\documentclass[a4paper,11pt]{article} 

\usepackage[utf8]{inputenc} 
\usepackage{microtype} 
\usepackage[english]{babel} 
\usepackage{amsthm, amsmath, amssymb} 
\usepackage{float} 
\usepackage[final, colorlinks = true, 
linkcolor = black, 
citecolor = black]{hyperref} 
\usepackage{caption, subcaption}
\usepackage{graphicx, multicol} 
\usepackage{xcolor} 
\usepackage{rotating} 
\usepackage{censor} 
\usepackage{matlab-prettifier} 
\usepackage{booktabs} 
\usepackage[yyyymmdd]{datetime} 
\usepackage{fancyhdr} 
\usepackage{pgfplots}
\pgfplotsset{compat=1.15}
\usepackage{hyperref}
\usepackage{tikz-cd,mathtools}

\usetikzlibrary{arrows}
\usetikzlibrary{braids}

\usepackage{wrapfig}
\usepackage{pdfpages}
\usepackage{mwe}

\usepackage{caption}

\DeclareMathOperator{\R}{\mathbb R}

\DeclareMathOperator{\N}{\mathbb N}
\DeclareMathOperator{\Q}{\mathbb Q}
\DeclareMathOperator{\Z}{\mathbb Z}

\DeclareMathOperator{\psl}{\mathbb P SL_2}
\DeclareMathOperator{\SL}{SL_2}

\DeclareMathOperator{\HH}{\mathbb H}

\DeclareMathOperator{\lk}{lk}
\DeclareMathOperator{\Ss}{S^2_*}
\DeclareMathOperator{\rays}{Rays}

\theoremstyle{plain}
\newtheorem{theorem}{Theorem}
\newtheorem{lemma}[theorem]{Lemma}
\newtheorem{proposition}[theorem]{Proposition}
\newtheorem{corollary}[theorem]{Corollary}
\newtheorem{conjecture}[theorem]{Conjecture}

\newtheorem{question}[theorem]{Question}
\newtheorem{fact}[theorem]{Fact}
\theoremstyle{definition}
\newtheorem{remark}[theorem]{Remark}

\newtheorem{assumptions}[theorem]{Assumptions}
\newtheorem{definition}[theorem]{Definition}
\newtheorem{notazione}[theorem]{Notation}
\newtheorem{example}[theorem]{Example}
\newtheorem{exercise}[]{Exercise}
\newtheorem{problem}[theorem]{Problem}

\newtheorem{vuoto}[theorem]{}

\numberwithin{theorem}{section}

\newcommand{\bt}{\begin{theorem}}
\newcommand{\et}{\end{theorem}}

\newcommand{\bv}{\begin{vuoto}}
\newcommand{\ev}{\end{vuoto}}

\newcommand{\bl}{\begin{lemma}}
\newcommand{\el}{\end{lemma}}

\newcommand{\bd}{\begin{definition}}
\newcommand{\ed}{\end{definition}}

\newcommand{\beq}{\begin{equation}}
\newcommand{\eeq}{\end{equation}}

\newcommand{\bexa}{\begin{example}}
\newcommand{\eexa}{\end{example}}

\newcommand{\bexe}{\begin{exercise}}
\newcommand{\eexe}{\end{exercise}}

\newcommand{\bfact}{\begin{fact}}
\newcommand{\efact}{\end{fact}}

\newcommand{\bprop}{\begin{proposition}}
\newcommand{\eprop}{\end{proposition}}

\newcommand{\bp}{\begin{proof}}
\newcommand{\ep}{\end{proof}}

\newcommand{\bc}{\begin{corollary}}
\newcommand{\ec}{\end{corollary}}

\newcommand{\bq}{\begin{question}}
\newcommand{\eq}{\end{question}}

\newcommand{\bcong}{\begin{conjecture}}
\newcommand{\econg}{\end{conjecture}}

\newcommand{\bproblem}{\begin{problem}}
\newcommand{\eproblem}{\end{problem}}

\newcommand{\bs}{\begin{proof}[Proof.]}
\newcommand{\es}{\end{proof}}

\newcommand{\br}{\begin{remark}}
\newcommand{\er}{\end{remark}}

\newcommand{\bn}{\begin{notazione}}
\newcommand{\en}{\end{notazione}}

\begin{document}

\title{The Thurston norm of 2-bridge link complements}
\author{Alessandro V. Cigna}
\date{\today}
\maketitle

\begin{abstract}
    \noindent The Thurston norm is a seminorm on the second real homology group of a compact orientable 3-manifold. The unit ball of this norm is a convex polyhedron, whose shape's data (e.g. number of vertices, regularity) measures the complexity of the surfaces sitting in the ambient 3-manifold. 

\noindent Unfortunately, the Thurston norm is generally quite hard to compute, and a long-standing problem is to understand which polyhedra are realised as the unit balls of the Thurston norms of $3$-manifolds. We show that, when $M$ is the complement of a $2$-bridge link $L$ with components $\ell_1$ and $\ell_2$, the Thurston ball of $M$ has at most 8 faces. The proof of this result strongly relies on a description of essential surfaces in $2$-bridge link complements given by Floyd and Hatcher in \cite{fh}. Then, we exhibit norm-minimizing representatives for the integral classes of $H_2(M,\partial M)$ and use them to compare the complexity of the Thurston ball with the complexities of $L$ and of $M$. As an example, we show that all the vertices of the Thurston ball lie on the bisectors if and only if $M$ fibers over the circle with fiber a surface with boundary equal to a longitude of $\ell_1$ and some meridians of $\ell_2$. Finally, we use $2$-bridge links in satellite constructions to find $2$-component links whose complements in $S^3$ have Thurston balls with arbitrarily many vertices. 
\end{abstract}

\section{Introduction}
The Thurston norm \cite{norm} is a powerful tool for analyzing the topology of a 3-manifold by studying the embedded surfaces sitting inside the $3$-manifold.
Given a compact orientable 3-manifold $M$, the Thurston norm is a seminorm $x$ defined on the vector space $H_2(M,\partial M;\R)$, whose value on a class represented by some properly embedded oriented surface $S$ gives the ``optimal" Euler characteristic among the properly embedded oriented surfaces homologous to $S$. For instance, in a knot complement the Thurston norm almost by definition detects the minimal genus of the knot. 

When $M$ is irreducible and $\partial M$ is a (possibly empty) union of tori, classical works by Thurston and Gabai \cite{norm, gabai1} establish that a properly embedded oriented surface $S\subset M$ with coherently-oriented boundary realizes the norm of its class if and only if $S$ is a union of leaves of a taut foliation of $M$. Here \emph{taut} means that, for every leaf, there is a closed transversal to the foliation intersecting that leaf. Taut foliations can then be used to study other surfaces embedded in $M$, by looking at the induced singular one-dimensional foliations on them \cite{norm, roussarie}.

The unit ball $B_x$ of $x$ is a finite (possibly unbounded) convex polyhedron (cf. \cite{norm}), thus the whole information carried by $x$ is synthesized by a finite amount of data, namely the position of the vertices of $B_x$. Moreover, the classes represented by fibers of fibrations of $M$ over $S^1$ lie in a union of open cones over top-dimensional faces of $B_x$ (cf. \cite{norm}). This and the subsequently developed theory of sutured manifolds (cf. \cite{gabai1}) were fundamental tools for solving conjectures that had remained open for a long time (e.g. \cite{agol, gabai1, gabai2, gabai3}). 

\bigskip

Unfortunately, the Thurston norm is generally quite hard to compute, and the Thurston ball has been described just in some special families of manifolds (see for example \cite{norm, pretzels, os, tunnel, chen}).  

If $L\subset S^3$ is an oriented link with components $\ell_1,...,\ell_n$, and $M=S^3-\overset{\circ}{N}(L)$ is the link exterior, then there is an isomorphism $H_2(M,\partial M;\Z)\to H_1(L;\Z)=(\Z\ell_1)\oplus...\oplus (\Z\ell_n)$ given by the map that associates to a surface $S$ the vector with $i$-th coordinate equal to the algebraic intersection number of $\partial S$ with the meridian of $\ell_i$.

In this article, we first compute the shape of the Thurston ball of the 3-manifolds $M$ obtained as complement of $2$-component $2$-bridge links $L_{p/q}$ in $S^3$, defined in Section \ref{subsec: links}. The main result of the paper is

\bigskip

\noindent\textbf{Theorem \ref{thm: Thurston ball}.}
\textit{Let $L_{p/q}$ be an oriented $2$-bridge link with components $\ell_1$ and $\ell_2$, and let $M=S^3-\overset{\circ}N(L_{p/q})$ be the link exterior.
   The Thurston unit ball in $H_2(M,\partial M;\R)=(\R\ell_1)\oplus (\R\ell_2)$ is the polygon spanned by the points $\pm\frac{\ell_1}{x(\ell_1)}$, $\pm\frac{\ell_2}{x(\ell_2)}$, $\pm\frac{\ell_1+\ell_2}{x(\ell_1+\ell_2)}$, and $\pm\frac{\ell_1-\ell_2}{x(\ell_1-\ell_2)}$.}

\bigskip

Basically, the vertices of $B_x$ are aligned with the somewhat most natural nontrivial classes in $H_2(M,\partial M;\R)$: the classes $\ell_1\pm\ell_2$ are represented by Seifert surfaces for the two (not projectively-equivalent) types of orientations of $L_{p/q}$, whilst the classes $\ell_1$ and $\ell_2$ are represented by Seifert surfaces for just one component, pierced by the other component in meridional discs.

Notice that the $2$-variable Alexander polynomial of $2$-bridge links has been computed in general (see \cite{hoste} for example). By \cite{os}, the Alexander polynomial determines the Thurston norm in alternating link complements. However, the arithmetic nature of the terms appearing in $2$-bridge link Alexander polynomials seems to obstruct a straightforward way to prove Theorem \ref{thm: Thurston ball} through these means. 

\bigskip

In order to compute the norms $x(\ell_1), x(\ell_2)$ and $x(\ell_1\pm \ell_2)$, we define properly embedded oriented surfaces $S_{1,0}$, $S_{0,1}$ and $S_{1,\pm 1}$ in $M=S^3-\overset{\circ}N(L_{p/q})$, representing respectively the classes $\ell_1, \ell_2$ and $\ell_1\pm \ell_2$. Then, for every $a,b\in \Z$, we define $S_{a,b}$ as the oriented cut-and-paste sum of some copies of the surfaces $S_{1,0}$, $S_{0,1}$ and $S_{1,\pm 1}$, possibly with inverted orientation (see Section~\ref{subsec: good representative}). 

\bigskip

\noindent\textbf{Theorem \ref{thm: standard representatives}.} \textit{For every $a,b\in \Z$, the surface $S_{a,b}(L)$ is a norm-minimizing representative for the class $a\ell_1+b\ell_2$.}

\bigskip

Of course, Theorem \ref{thm: Thurston ball} implies that the Thurston ball of a $2$-bridge link complement has at most $8$ sides. Whenever three of the vertices mentioned in Theorem \ref{thm: Thurston ball} happen to be aligned, a fewer number of vertices spans the Thurston ball, and we would expect such a lower complexity of the norm to come from a lower complexity of the link itself. This is indeed the case, as shown in Section \ref{subsec: simpler}. For instance, by Corollary \ref{cor: base-type}, $B_x$ has vertices along the bisectors if and only if $M$ fibers over the circle with fiber $S_{1,0}$.

\bigskip

A long-standing problem is to understand which polyhedra are realised as the unit balls of the Thurston norms of $3$-manifolds. In \cite{norm}, Thurston showed that every symmetric polygon with vertices in $\Z^2$ all entry-wise congruent modulo $2$ is the unit ball of the dual Thurston norm on $H^2(N,\partial N;\R)$, for some $3$-manifold $N$. Unfortunately, the family of manifolds considered by Thurston are not link complements in $S^3$, and the problem of understanding which polyhedra are realised as Thurston balls of link complements is still open. As shown by \cite{karim2}, the construction of Thurston can be generalized to manifolds with higher homological rank, but then many polyhedra satisfying the same minimal assumptions as before cannot be realized through the same construction anymore (cf. \cite{karim1}).

However, also in the $2$-dimensional homology setting, it is hard to keep track of the complexity of the manifolds $N$. In \cite{tunnel}, polygons with arbitrarily many faces are proved to be realized as Thurston balls of manifolds with two generators and one relator fundamental group.

In Section \ref{sec: satellites}, we use an iterated satellite construction to show a similar result, for link exteriors $M_L$ in $S^3$.

\bigskip

\noindent\textbf{Theorem \ref{thm: complexity}.} \textit{For every $n\in\N$, there is a $2$-component link $L\subset S^3$ with each component being the unknot, such that the Thurston ball of $H_2(M_L,\partial M_L;\R)$ has exactly $2n$ faces.}




\bigskip

\noindent \textbf{Structure of the paper.}  In Section \ref{sec: preliminaries} we recall the definitions and some of the main properties related to the Thurston norm, sutured manifolds and disc decompositions, and to $2$-bridge links. In Section \ref{sec: thurston ball} we slightly sharpen the work \cite{fh} by Floyd and Hatcher to prove Theorem \ref{thm: Thurston ball}. In Section \ref{sec: norm-minimizing} we first describe norm-minimizing surfaces in $2$-bridge link complements, then use them to compare the complexity of the Thurston ball with the complexity of the link. As an application of the previous sections, in Section \ref{sec: satellites}, we relate the norm of a satellite exterior to the norm of its summands and use this construction to find $2$-component links with arbitrarily complicated Thurston ball.

\bigskip

\noindent \textbf{Aknowledgements.} I would like to thank my advisor Mehdi Yazdi for suggesting that $2$-bridge link exteriors have Thurston ball with at most eight faces and for the enlightening discussions on the topic. I would also like to thank Steven Sivek for his useful advice during the preparation of this work.

\section{Preliminaries}\label{sec: preliminaries}
\subsection{Thurston norm}
We now recall the basic definitions and properties of the Thurston norm. For further details, we refer to \cite{norm}. From now on, $M$ will always be a $3$-dimensional compact orientable manifold.

\bigskip

\noindent Every properly embedded oriented surface $S\subset M$ represents a class $[S]\in H_2(M,\partial M; \Z)$. Conversely, every class in $H_2(M,\partial M;\Z)$ is represented by a properly embedded oriented surface. 

It's quite a natural question to ask whether the sum in $H_2(M,\partial M;\Z)$ corresponds with some operation between surfaces in $M$. The answer is affirmative: if $[S_1]=\alpha_1$ and $[S_2]=\alpha_2$, then $\alpha_1+\alpha_2$ is represented by the \emph{oriented cut-and-paste} of $S_1$ and $S_2$. Here is the definition.

\bd(Oriented cut-and-paste) Let $S_1,S_2\subset M$ be properly embedded oriented surfaces. The \emph{oriented cut-and-paste} (or the \emph{double-curve sum}) of $S_1$ and $S_2$ is (up to isotopy) the surface $S$ defined as follows.

\begin{figure}[ht]
    \centering
    \includegraphics[width=0.85\textwidth]{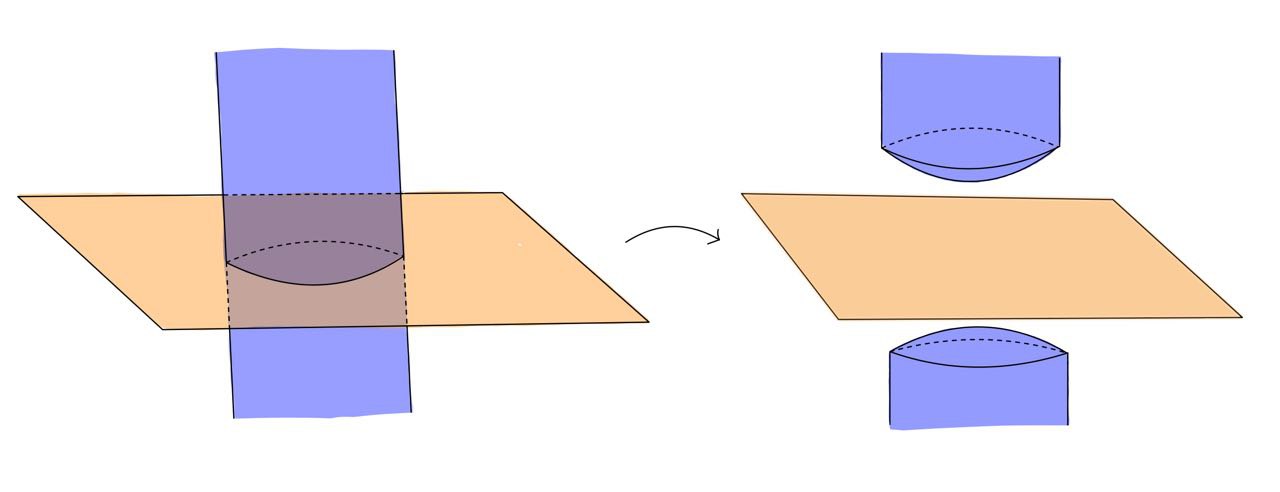}
    \caption{Surgery along an inessential intersection. In blue the surface $S_1$, in orange $S_2$.}
    \label{fig: 2-surgery}
\end{figure}

\begin{figure}[ht]
    \centering
    \includegraphics[width=0.85\textwidth]{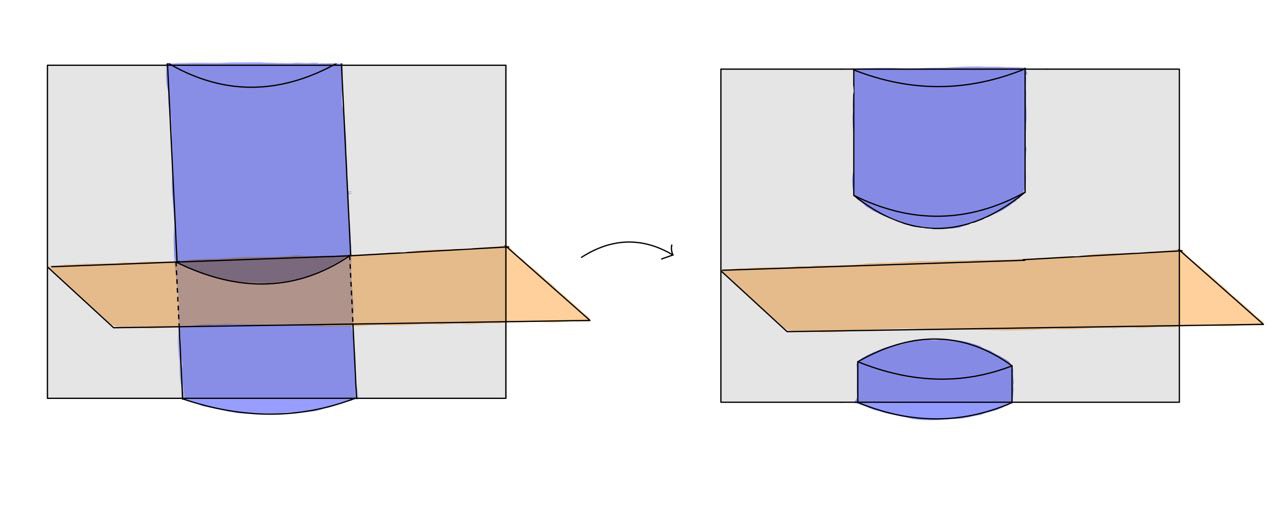}
    \caption{Boundary surgery along an inessential intersection. In blue and orange the surfaces $S_1$ and $S_2$, in grey $\partial M$.}
    \label{fig: boundary surgery}
\end{figure}

\begin{figure}[ht]
    \centering
    \includegraphics[width=0.85\textwidth]{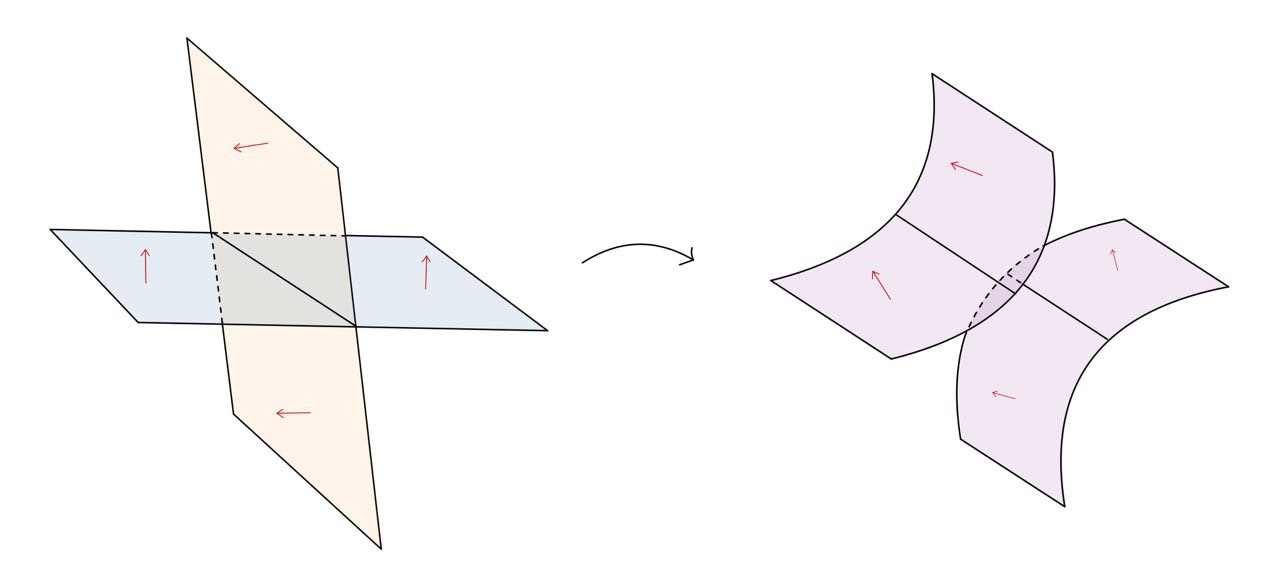}
    \caption{Local view of the oriented cut-and-paste operation of the surfaces $S_1$ and $S_2$. Red arrows indicate the normal orientation.}
    \label{fig: orientedcutandpaste}
\end{figure}

If $S_1\cap S_2$ consists of some circle which bounds a disc $D$ in $S_1$, then do a surgery on $S_1$ along an innermost circle of $S_1\cap S_2$ contained in $D$, see Figure \ref{fig: 2-surgery}. By iterating this procedure, also exchanging the roles of $S_1$ and $S_2$ if needed, we modify $S_1$ and $S_2$ up to assume that no circle in $S_1\cap S_2$ bounds a disc in either surface. Using a similar technique, we can also assume that no properly embedded arc in $S_1\cap S_2$ is boundary-parallel (see Figure \ref{fig: boundary surgery}). Finally, there is only one oriented way to split $S_1$ and $S_2$ along $S_1\cap S_2$ and reglue back to obtain the oriented surface $S$, which locally looks as in Figure \ref{fig: orientedcutandpaste}.
\ed

Notice that the standard definition of oriented cut-and-paste does not require to eliminate the essential intersections between $S_1$ and $S_2$ first. We incorporate this in the definition so that the triangular inequality for the Thurston norm can be verified using $S_1, S_2$ and their oriented cut-and-paste.

\bd(Thurston-norm) Let $S\subset M$ be a properly embedded orientable surface. Denote with $\chi(S)$ the Euler characteristic of $S$. If $S$ is connected, call $$\chi_-(S):=\max (0,-\chi(S)).$$
If $S$ is the disjoint union of connected surfaces $S_1,...,S_k$, define $$\chi_-(S):=\chi_-(S_1)+...+\chi_-(S_k).$$
For every class $\alpha\in H_2(M,\partial M;\Z)$, let $$x(\alpha)=\min_{[S]=\alpha}\chi_-(S),$$ where $S\subset M$ varies among the properly embedded oriented surfaces representing $\alpha$.
\ed

\bt[\cite{norm}] The function $x:H_2(M,\partial M;\Z)\to \Z$ can be extended to a seminorm $x:H_2(M,\partial M;\R)\to \R_{\ge 0}$, called the \emph{Thurston norm} associated to the manifold $M$.
\et
\noindent\emph{Sketch of the proof.} Thurston verifies that $x$ is $\N$-homogeneous on $H_2(M,\partial M;\Z)$ and that $x(S)\le x(S_1)+x(S_2)$ whenever $S$ is the oriented cut-and-paste of $S_1$ and $S_2$. Then, $x$ can be extended by homogeneity to a convex function $x: H_2(M,\partial M;\Q)\to \Q$, which can be continuously extended further to a seminorm $x: H_2(M,\partial M;\R)\to \R_{\ge 0}$. 

\bigskip

To compute a (semi)norm on $\R^n$ it is enough to know its unit ball. A handy feature of Thurston norm is that its ball can be recovered by knowing a finite amount of data, thanks to the following result.

\bt[\cite{norm}] The unit ball $B_x$ of the Thurston norm $x$ on $M$ is a finite (but possibly unbounded) convex polyhedron, symmetric with respect to the origin, and whose vertices are rational points in $H_2(M,\partial M;\R)$. 
\et

When $x$ is not a norm, i.e. when some surface with non-negative Euler characteristic is not null-homologous, the null space of $x$ is a vector space contained in $B_x$. In this case, in the previous theorem we allow vertices to be points at infinity along a ray from the origin, i.e. points in the visual boundary of $H_2(M,\partial M;\R)$. Of course, $B_x$ is bounded if and only if $x$ is a norm on $H_2(M,\partial M;\R)$. 

\bigskip

A useful principle to investigate the shape of the Thurston ball of a manifold is the following.

\bl\label{lemma: additivity} Let $\alpha,\beta\in H_2(M,\partial M;\R)$. Then $\alpha$ and $\beta$ lie in the cone over a common face of the Thurston ball if and only if $$x(\alpha+\beta)=x(\alpha)+x(\beta).$$
\el
\bp If $\alpha$ and $\beta$ lie in the cone over a common face $F$ of the Thurston ball, then every point of the segment between $\alpha/x(\alpha)$ and $\beta/x(\beta)$ lies on $F$. In particular \begin{equation}\label{eq: additivity} x\left(\frac{x(\alpha)}{x(\alpha)+x(\beta)} \frac{\alpha}{x(\alpha)}+\frac{x(\beta)}{x(\alpha)+x(\beta)}\frac{\beta}{x(\beta)}\right)=1,\end{equation} and this is equivalent to $x(\alpha+\beta)=x(\alpha)+x(\beta)$.

Viceversa, $x(\alpha+\beta)=x(\alpha)+x(\beta)$ implies the validity of equation (\ref{eq: additivity}). Thus, the three aligned points $$\frac{\alpha}{x(\alpha)},\frac{\beta}{x(\beta)} \text{ and } \frac{\alpha+\beta}{x(\alpha+\beta)}=\frac{x(\alpha)}{x(\alpha)+x(\beta)} \frac{\alpha}{x(\alpha)}+\frac{x(\beta)}{x(\alpha)+x(\beta)}\frac{\beta}{x(\beta)}$$ all lie on the boundary of the Thurston ball, so the whole segment between $\alpha/x(\alpha)$ and $\beta/x(\beta)$ does.  
\ep

\subsection{Sutured manifolds and fibering classes}
Sutured manifolds were introduced by Gabai \cite{gabai} to construct taut foliations on $3$-manifolds. Indeed, by these means, he was able to show that every compact orientable irreducible $3$-manifold $M$ with toroidal boundary admits taut foliations whenever $M$ has non-null first Betti number \cite{gabai1}.

For our purposes, just one tool from the theory of sutured manifolds is needed, namely disc decompositions. However, thanks to Theorem \ref{thm: fibering}, this single tool is already enough to characterize when a $3$-manifold fibers over $S^1$ with a given fiber. 

\bd[Sutured manifold] A \emph{sutured manifold} is a pair $(M,\gamma)$, where $M$ is a compact oriented $3$-manifold and $\gamma=A(\gamma)\cup T(\gamma)$ is a subset of $\partial M$ satisfying the following conditions. The set $A(\gamma)$ is a union of pairwise disjoint annuli, whereas $T(\gamma)$ is a union of pairwise disjoint tori. A \emph{suture} is the core of an annulus of $A(\gamma)$, considered as an oriented simple closed curve. Every region of $R(\gamma):=\partial M-\overset{\circ}{\gamma}$ has a normal orientation that can disagree with the one of $\partial M$. We call $R_+(\gamma)$ the union of the components of $R(\gamma)$ whose normal vector points out of $M$, and $R_-(\gamma)$ the union of the remaining components. The last requirement for $(M,\gamma)$ is that the normal orientations on $R(\gamma)$ are coherent with the orientation of the sutures: if $\delta\subset \partial R(\gamma)$ is a boundary component with the induced orientation, then $\delta$ is homologous to a suture in $H_1(\gamma)$.
\ed

The idea behind sutured manifolds is simple: assume $S$ is a properly embedded oriented surface in a manifold $M$ with toroidal boundary, this means that $(M,\partial M)$ is a sutured manifold. Suppose that we want to construct a foliation of $M$ where $S$ is a leaf. Cut $M$ open along $S$, so to obtain the manifold $N=M-\overset{\circ}{N}(S)$. We can now consider the sutured manifold $(N,\gamma)$, with $\gamma=\partial M-\overset{\circ}{N}(S)$ and sutures isotopic to $\partial S$. Here, both $R_+(\gamma)$ and $R_-(\gamma)$ are copies of $S$. So, if we find a foliation of $N$ transverse to $\gamma$ and tangent to $R_+(\gamma)$ and $R_-(\gamma)$, then we can construct a foliation of $M$ with leaf $S$, by just gluing back $R_+(\gamma)$ and $R_-(\gamma)$.

Potentially, this technique can be reiterated: we find a new surface $S'$ in $N$ and cut $N$ open along $S'$. Then we define a new sutured manifold $(N',\gamma')$, where $N'=N-\overset{\circ}{N}(S')$ and $\gamma'$ is defined so that a foliation of $N'$ transverse to $\gamma'$ and tangent to $R(\gamma')$ can be transformed to get a foliation of $(N,\gamma)$ as above. Then we can decompose $(N',\gamma')$ along a new surface $S''$ and so on and so forth.

Finally, suppose that at some point we get a sutured manifold homeomorphic to $(\Sigma\times [0,1], \partial\Sigma\times [0,1])$, where $\Sigma$ is a surface with non-empty boundary. Then the process is done, because we can just choose the product foliation on this last sutured manifold.

This philosophy stands behind the notions of \emph{sutured manifold decompositions} and \emph{sutured manifold hierarchies}. As previously mentioned, we limit ourselves to product-discs decompositions.

\bd(Product-discs decompositions) Let $(M,\gamma)$ be a sutured manifold. A \emph{product-disc} for $(M,\gamma)$ is a properly embedded oriented disc $D\subset M$ that intersects transversely exactly two sutures of $\gamma$.

\begin{figure}[h]
    \centering
    \includegraphics[width=0.75\textwidth]{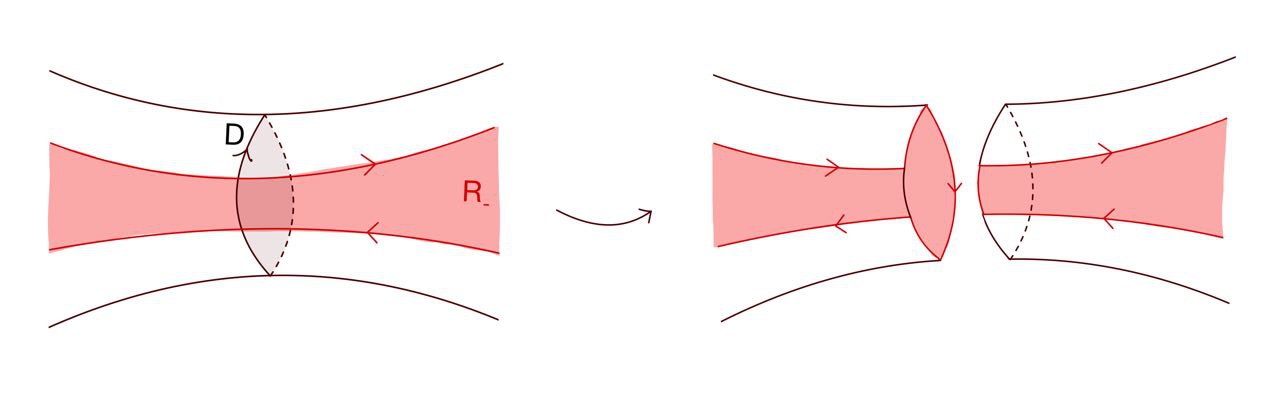}
    \caption{}\label{fig: decomposition}
\end{figure}

Given a product-disc $D\subset M$, we can \emph{decompose} $M$ along $D$ in order to obtain a new sutured manifold $(M',\gamma')$, where $M'=M-\overset{\circ}{N}(D)$ and $\gamma'$ is equal to $\gamma$ out of $\partial N(D)$, and is modified as in Figure \ref{fig: decomposition} on $\partial N(D)$. We indicate product-disc decompositions as $(M,\gamma)\rightsquigarrow (M',\gamma')$. 

We say that $(M,\gamma)$ is \emph{completely decomposable through product-discs}, if there is a sequence of product-disc decompositions $(M,\gamma)\rightsquigarrow (M',\gamma')\rightsquigarrow...\rightsquigarrow (M^{(k)},\gamma^{(k)})$, such that $(M^{(k)},\gamma^{(k)})$ is homeomorphic to a disjoint union of sutured balls $(D^2\times[0,1],\partial D^2\times [0,1])$.
\ed

\bt\cite{gabai}\label{thm: fibering} Let $(M,\gamma)$ be a sutured manifold with $R(\gamma)\ne \emptyset$. If $(M,\gamma)$ is completely decomposable through product-discs, then $(M,\gamma)$ is a product-sutured manifold: i.e. $(M,\gamma)$ is homeomorphic to $(R_+(\gamma)\times[0,1],\partial R_+(\gamma)\times [0,1])$.
\et

\subsection{$2$-bridge links}\label{subsec: links}
This subsection is dedicated to fixing the notations concerning $2$-bridge links. There are two principal ways to visualize these links: rational diagrams and $2$-bridge positions. If the latters allow us to understand the shape of the Thurston ball, the formers are the ones that make real computations of the Thurston norm possible.

\begin{definition}
    Let $a_1,...,a_k$ be non-zero integers. The \emph{rational diagram} $T(a_1,...,a_k)$ is the link diagram in Figure \ref{fig: rational} top and centre. Here, each box labeled $a\in\Z$ represents a sequence of $|a|$ crossings, as in Figure \ref{fig: rational} bottom. 

    \begin{figure}[ht]
    \centering
    \includegraphics[width=0.7\textwidth]{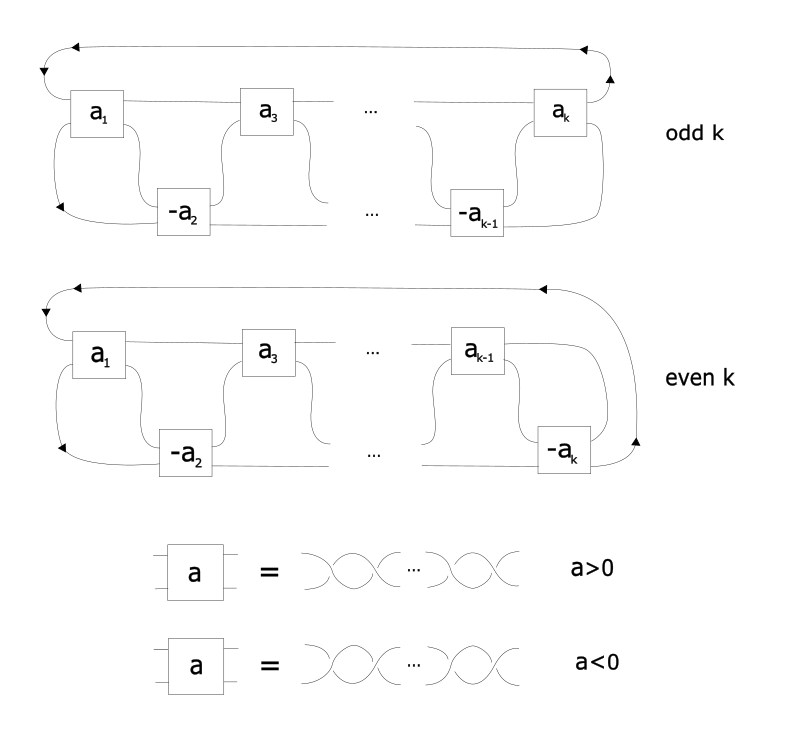}
    \caption{}\label{fig: rational}
\end{figure}

    \noindent Every rational diagram comes with a \emph{standard orientation}, i.e. the one represented in Figure \ref{fig: rational}.

    A \emph{$2$-bridge link} is a link that admits a rational diagram. 
\end{definition}

\noindent
Thanks to a classic work by Schubert \cite{schubert} (see also \cite{bz}), two non-oriented rational diagrams $T(a_1,...,a_k)$ and $T(b_1,...,b_h)$ represent the same 2-bridge link if the associated continued fraction expansions are equal, i.e. if $$[a_1,...,a_k]:=\frac{1}{a_1+\frac{1}{a_2+\frac{1}{...\frac 1{a_k}}}} = \frac{1}{b_1+\frac{1}{b_2+\frac{1}{...\frac 1{b_h}}}} =:[b_1,...,b_h].$$ Given a rational number $p/q\in (0,1)$, we call $L_{p/q}$ the $2$-bridge link admitting a rational diagram with associated continued fraction $p/q$.

There is an isotopic representative for $L_{p/q}$ depending only on the fraction $p/q$. In order to describe it, consider the group $\Gamma$ generated by 180-degree rotations around $\Z^2$ points in $\R^2$. The quotient $\R^2/\Gamma$ is homeomorphic to a sphere $S^2$. We will refer to the surface $\R^2/\Gamma$ with the four marked points $\Z^2/\Gamma$ as $S^2_*$. Every affine line in $\R^2$ with rational slope and passing through a point of $\Z^2$ is mapped to a simple arc in $S^2_*$ with two distinct marked points as endpoints. In fact, given a marked point $x\in\Ss$, any simple arc with endpoints $x$ and a different marked point is isotopic (rel. its endpoints) to exactly one such arc with rational slope.

\begin{figure}[ht]
    \centering
    \includegraphics[width=0.5\textwidth]{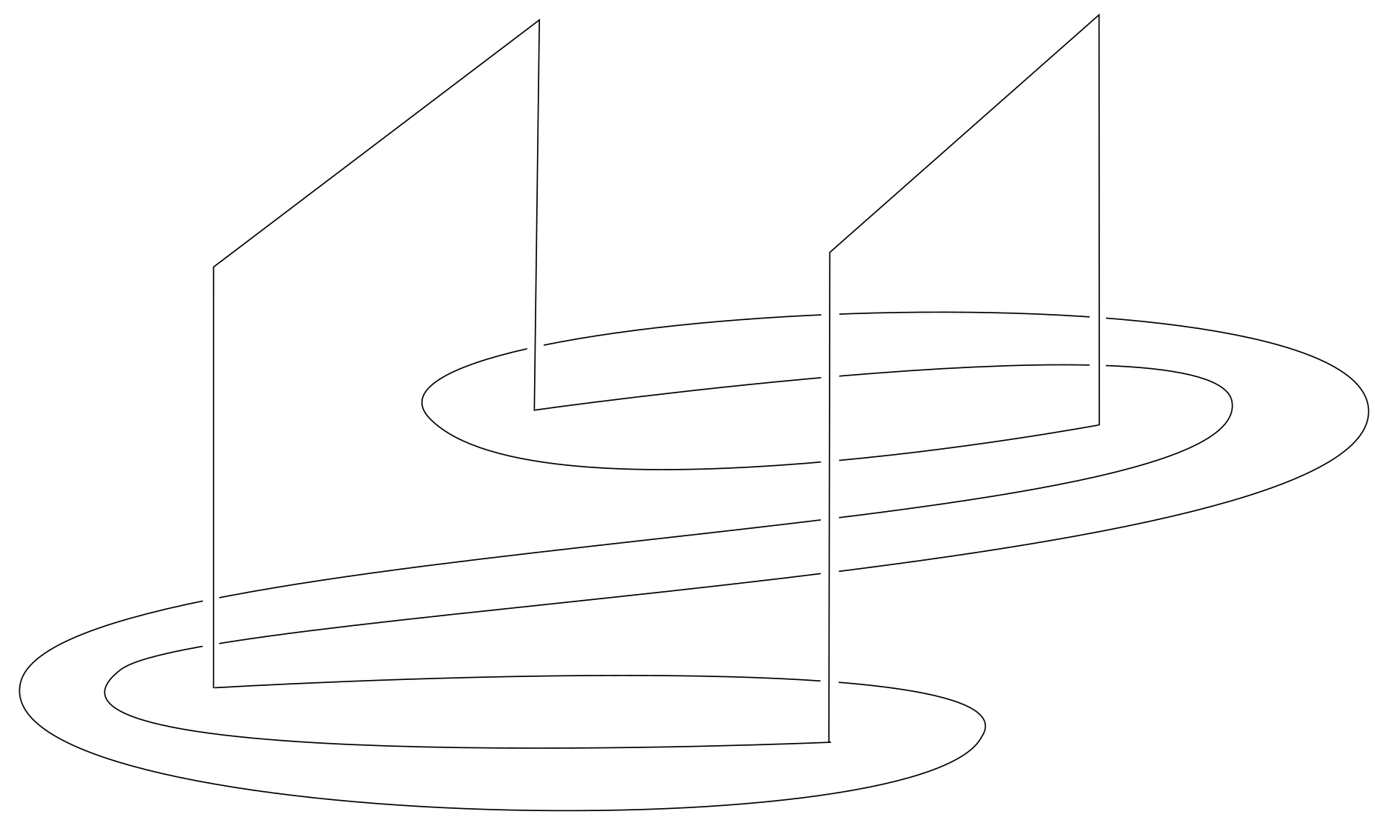}
    \caption{Schubert normal form for $L_{1/4}$.}\label{fig: schubert normal form}
\end{figure}

The \emph{$2$-bridge position} for $L_{p/q}$ is the link $L\subset \Ss\times [0,1]\subset S^3$ such that \begin{itemize}
    \item $L\cap \Ss\times\{1\}$ is the union of two disjoint arcs of slope $1/0$ each connecting a different pair of marked points;
    \item $L\cap \Ss\times\{t\}$ consists exactly of the four marked points, for every $t\in (0,1)$;
    \item $L\cap \Ss\times\{0\}$ is the union of two disjoint arcs of slope $p/q$ each connecting a different pair of marked points (see Figure~\ref{fig: schubert normal form}).
\end{itemize}

From Schubert's work, any $2$-bridge link is alternating (by choosing a continued fraction expansion with all $a_i$ of the same sign), and $L_{p/q}$ is a knot if and only if $q$ is odd, whereas $L_{p/q}$ has exactly two components for even $q$.
We will only work with two-component $2$-bridge links, i.e. with even $q$, that we will refer to just as $2$-bridge links.

\section{The Thurston ball of $2$-bridge link complements}\label{sec: thurston ball}
In this section, we describe the Thurston ball associated with $2$-bridge link complements. In particular, this happens to be a (possibly degenerate) octagon. A cone over one of its edges is a union of some of the 8 cones given by cutting $\R^2$ along the axes and the principal bisectors. The proof is achieved essentially via Morse theory, after finding a nice norm-minimizing representative for each integral homology class.

We will consistently use the work by Floyd and Hatcher \cite{fh}, where they classify the properly embedded incompressible, orientable surfaces in the link-complement $M_{L_{p/q}}=S^3-\overset{\circ}{N}(L_{p/q})$ up to isotopy. 

For the rest of the section, let $L_{p/q}\subset \Ss\times [0,1]$ be in $2$-bridge position and call $M=M_{L_{p/q}}$. Observe that since $L_{p/q}$ is nontrivial and nonsplit, then $M$ is irreducible and $\partial$-irreducible, and so $\chi_-=-\chi$ on every connected surface which is not null-homologous.

\subsection{A good representative for each homology class}\label{subsec: good representative}

\begin{assumptions}\label{assumptions}
Let $S\subset M$ be a properly embedded oriented surface. Through the identification $H_2(M,\partial M)\cong H_1(L_{p/q})$, we may orient and label the link components $\ell_1$ and $\ell_2$ such that $[S]=a\ell_1+b\ell_2$ with $0\ne a\ge b\ge 0$. We now make some modifications on $S$ to track better how it is embedded in $M$.

\begin{itemize}
    \item No null-homologous components in $\partial S$. After possibly sewing in discs to $\partial S$, we can assume that no component of $\partial S$ is null-homologous in the component of $\partial M$ where it lies.
    \item Meridional-incompressibility. Suppose there is an annulus $A\subset M$ with $\partial A=A\cap (S\cup\partial M)$, one component of $\partial A$ lying on $S$ and one isotopic to a meridian of $\partial M$. In that case, we can operate a \emph{meridional-surgery} as in Figure~\ref{fig: compression}. After a finite number of surgeries, we can suppose $S$ to be meridionally-incompressible, meaning that any annulus $A$ as above is parallel (rel. $\partial A$) to one contained in $S$.
    \item Coherent orientation of $\partial S$ apart for meridians. After possibly sewing in some annuli to $\partial S$ and isotopy, we can assume that the simple closed curves of $S\cap \partial N(\ell_1)$ represent the same homology class and that they are transverse to every meridian. If $b\ne 0$, we make the same assumption on $S\cap\partial N(\ell_2)$. If $b=0$, then $S\cap\partial N(\ell_2)$ consists of meridians, and the annular surgery would be incompatible with the meridional-incompressibility. In this case, we do not ask the meridians to be coherently oriented.
    \item Incompressibility. After possibly surgeries along compression discs, we can assume $S$ is incompressible.
    \item No closed components. Every closed surface in $M$ is null-homologous in $H_2(M,\partial M)$, so we can discard any closed component from $S$.
    \item No peripheral components. If a component of $S$ is $\partial$-parallel, then it is homologically trivial, so we discard it from $S$.
\end{itemize}

\begin{figure}[ht]
    \centering
    \includegraphics[width=\textwidth]{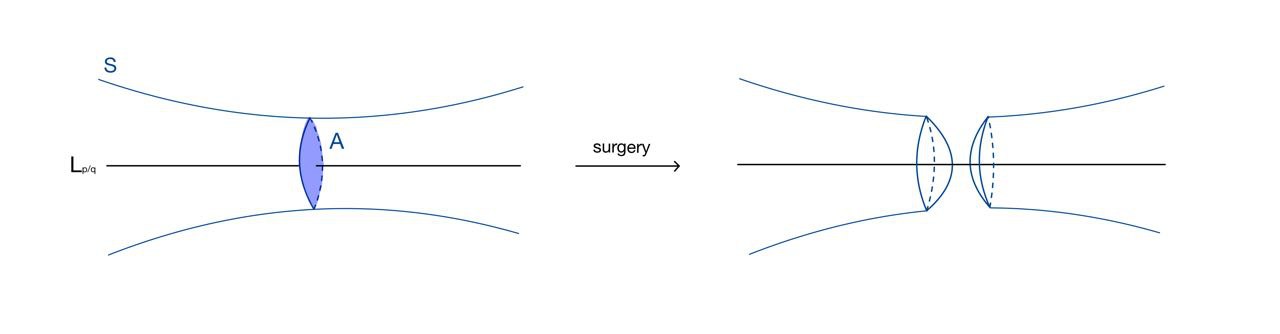}
    \caption{}\label{fig: compression}
\end{figure}

Please notice that none of the above transformations increases $\chi_-(S)$ and, in particular, we can find a norm-minimizing representative for $a\ell_1+b\ell_2$ satisfying all the conditions above. Also, the transformations above should be done in that order so as not to interfere ones with each other (e.g. compressions can create closed components). 
\end{assumptions}
 
\bigskip

Since $L_{p/q}$ is in $2$-bridge position, the height function $h:\Ss\times [0,1]\to [0,1]$ gives a good description of the sections of $L_{p/q}$ through the horizontal spheres $\Ss\times\{t\}$. We would like a similar description of the sections of $S$ through the same horizontal spheres. For this purpose, we will conveniently think of the tubular neighborhood $N(L_{p/q})$ to be as narrow as to identify $\overset{\circ}{M}$ with $S^3-L_{p/q}$, and $S$ to be a compact singular surface in $S^3$ with $\partial S\subset L_{p/q}$ and embedded interior. Thanks to the coherent orientation of $\partial S$, we know that $\partial S$ wraps $a$ times along $\ell_1$ and, if $b\ne 0$, $b$ times along $\ell_2$. When $b=0$, the surface $S$ intersects $\ell_2$ transversely. 

We can now apply the next result, implicitly shown in the proof of Theorem $3.1$ of \cite{fh}.

\bt\cite[Section 7]{fh}\label{thm: fh1} Let $S\subset M$ be a properly embedded oriented surface. If $S$ is incompressible, meridionally-incompressible, and has no peripheral components, then we can isotope $S$ (rel. its boundary) so that it is contained in $\Ss\times [0,1]$, $S\cap \Ss\times \{0,1\}\subset L_{p/q}$, and the height function $h:\overset{\circ}{S}\cap\Ss\times(0,1)\to (0,1)$ is a generic Morse function with only saddle singularities.
\et

For every height $t\in [0,1]$ call $S_t=S\cap\Ss\times \{t\}$. For any non-singular height $t$, the level set $S_t\subset \Ss$ is a disjoint union of some simple closed curves and an arc-system. The arc-system is simply a graph whose vertices are the marked points of $\Ss$. Our assumption on the orientation of $\partial S$ implies that two vertices of $\Ss$ have always degree $a$, and the other two have always degree $b$. 

\br Fixing a partition of the marked points of $\Ss$ into two pairs, call $\Gamma(a,b)$ the set of possible graphs in $\Ss$ with vertices in the marked points, the first couple of them with degree $a$ and the second with degree $b$.

The fundamental reason why $x(2\ell_1+\ell_2)=x(\ell_1)+x(\ell_1+\ell_2)$ - and so the segment with endpoints $\ell_1/x(\ell_1)$ and $(\ell_1+\ell_2)/x(\ell_1+\ell_2)$ lies on the boundary of the Thurston ball - is that every arc-system in $\Gamma(2,1)$ can be split as the disjoint union of two arc-systems, one lying in $\Gamma(1,0)$ and the other in $\Gamma(1,1)$. Indeed, Theorem \ref{thm: Thurston ball} can be shown without introducing the following machinery (the complexes $D_\alpha$ and their subtrees $T_{\alpha}$). This alternative proof can be obtained by a reasoning similar to the one in Remark \ref{rem: cut-and-paste}, where we show that a nice norm-minimizing representative for $2\ell_1+\ell_2$ must be the oriented cut-and-paste of some representatives for $\ell_1$ and $\ell_1+\ell_2$. However, what follows will be needed in the proof of Theorem \ref{thm: standard representatives} as well.
\er

\begin{figure}[ht]
    \centering
    \includegraphics[width=0.3\textwidth]{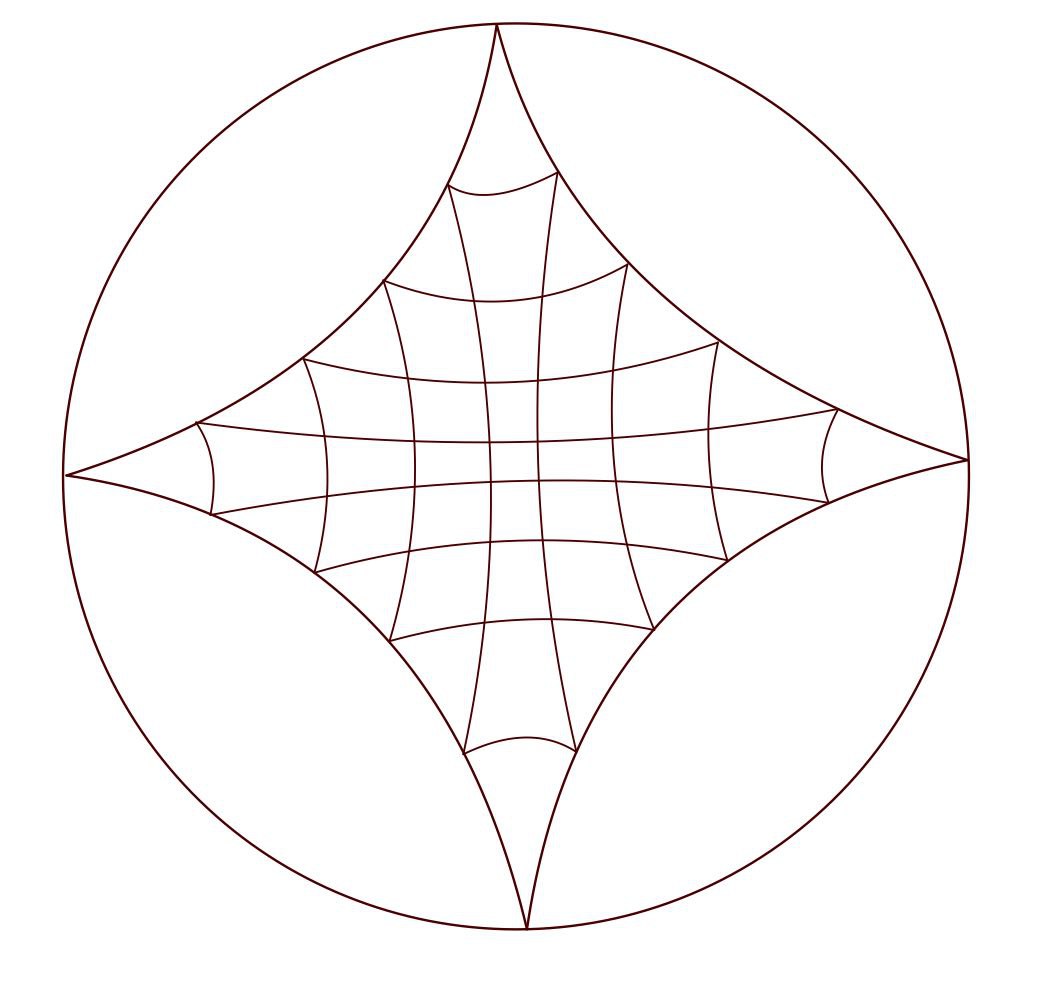}
    \caption{The inscribed quadrilateral changes its shape varying with $\alpha$. For $\alpha\to 0$ or $\infty$ it tends to the horizontal diagonal, for $\alpha\to 1$ it tends to the vertical one.}\label{fig: quadrilaterals}
\end{figure}

\subsection{The complexes $D_{\alpha}$}

\begin{figure}[hb]
    \centering
    \includegraphics[width=0.7\textwidth]{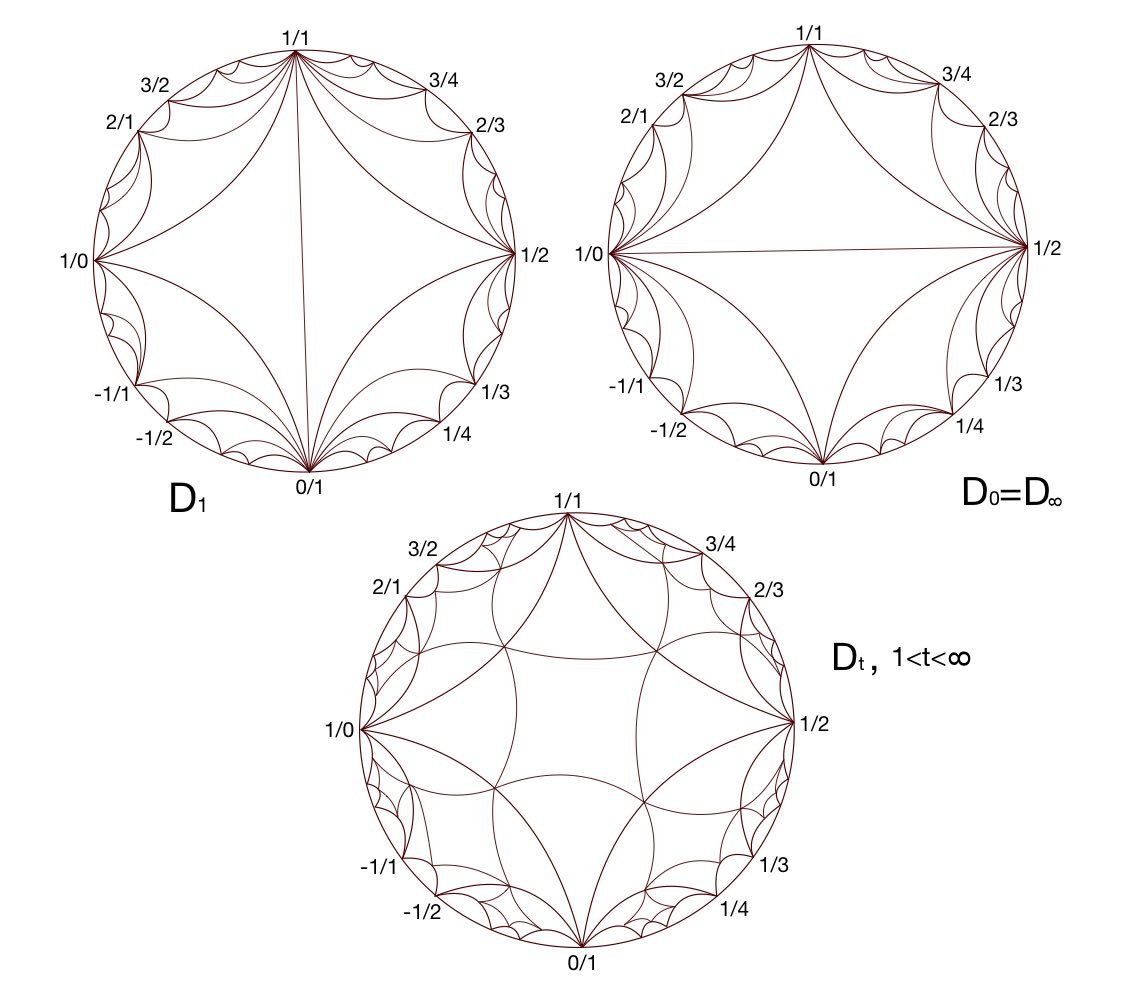}
    \caption{}\label{fig: diagrams}
\end{figure}

Floyd and Hatcher in \cite{fh} exhibit a $2$-complex $D_\alpha$, for every rational $\alpha\in \Q\cup\{\infty\}$, with the idea that the finite sequence of non-isotopic sections $S_t$ should indicate a path on the $1$-skeleton of $D_\alpha$, for some $\alpha$. In our case, we will have $\alpha=a/b$, but let us first introduce the complexes $D_{\alpha}$, shown in Figure \ref{fig: diagrams}.

\bigskip

In this subsection, we almost textually cite the construction of \cite[section $1$]{fh}.
The diagram $D_1$ is the Farey tessellation of the hyperbolic plane. The vertices of $D_1$ are the points $\Q\cup\{\infty\}\subset \partial \HH^2$, an edge joins two reduced fractions $c/d, e/f$ if and only if $|cf-de|=1$. The group $\psl (\Z)$ is the full group of orientation-preserving symmetries of this ideal triangulation. 
Let $G \subset \psl(Z)$ be the subgroup of M\"obius transformations $(cz+d)/(ez+f)$ with $e$ even. This has index three in $\psl (Z)$ and has the ideal triangle $(1/0, 0/1, 1/1)$ as a fundamental domain. The element $z+ 1\in G$ identifies the edge $(1/0, 0/1)$ to the edge $(1/0,1/1)$, while $(z- 1)/(2z- 1)\in G$ takes
 the edge $(0/1, 1/1)$ to itself, reversing its endpoints. Consider the ideal quadrilateral $(1/0, 0/1, 1/2, 1/1)$, which is rotated $180$ degrees about its center by $(z - 1)/(2z - 1)$. The
 $G$-images of this quadrilateral form a tiling of $\HH^2$ by ideal quadrilaterals. We
 form the diagram $D_0$ from $D_1$ by deleting the $G$-orbit of the diagonal $(0/1, 1/1)$
 of $(1/0, 0/1,1/2,1/1)$ and adding the $G$-orbit of the opposite diagonal $(1/0,1/2)$.
 Between $D_0$ and $D_1$ we interpolate a $1$-parameter family of $G$-invariant diagrams
 $D_{\alpha}$, $0 < \alpha < 1$, by inscribing a rectangle in $(1/0,0/1,1/2,1/1)$ of monotonically
 varying shape, as in Figure \ref{fig: quadrilaterals}, the rectangle collapsing to the diagonals $(1/0,1/2)$
 and $(0/1,1/1)$ as $\alpha$ approaches $0$ and $1$, respectively. Such inscribed rectangles are
 determined by their vertices, which we take to be the $G$-orbit of a point in the edge
 $(1/0, 0/1)$, this point moving monotonically from 
 $1/0$ to $0/1$ as $\alpha$ goes from $0$ to $1$.
 This defines $D_{\alpha}$ for $\alpha\in [0,1]$, and we obtain $D_\alpha$ for $\alpha\in [1, \infty]$ by setting $D_\alpha = D_{1/\alpha}$.

\begin{figure}[ht]
    \centering
    \includegraphics[width=0.8\textwidth]{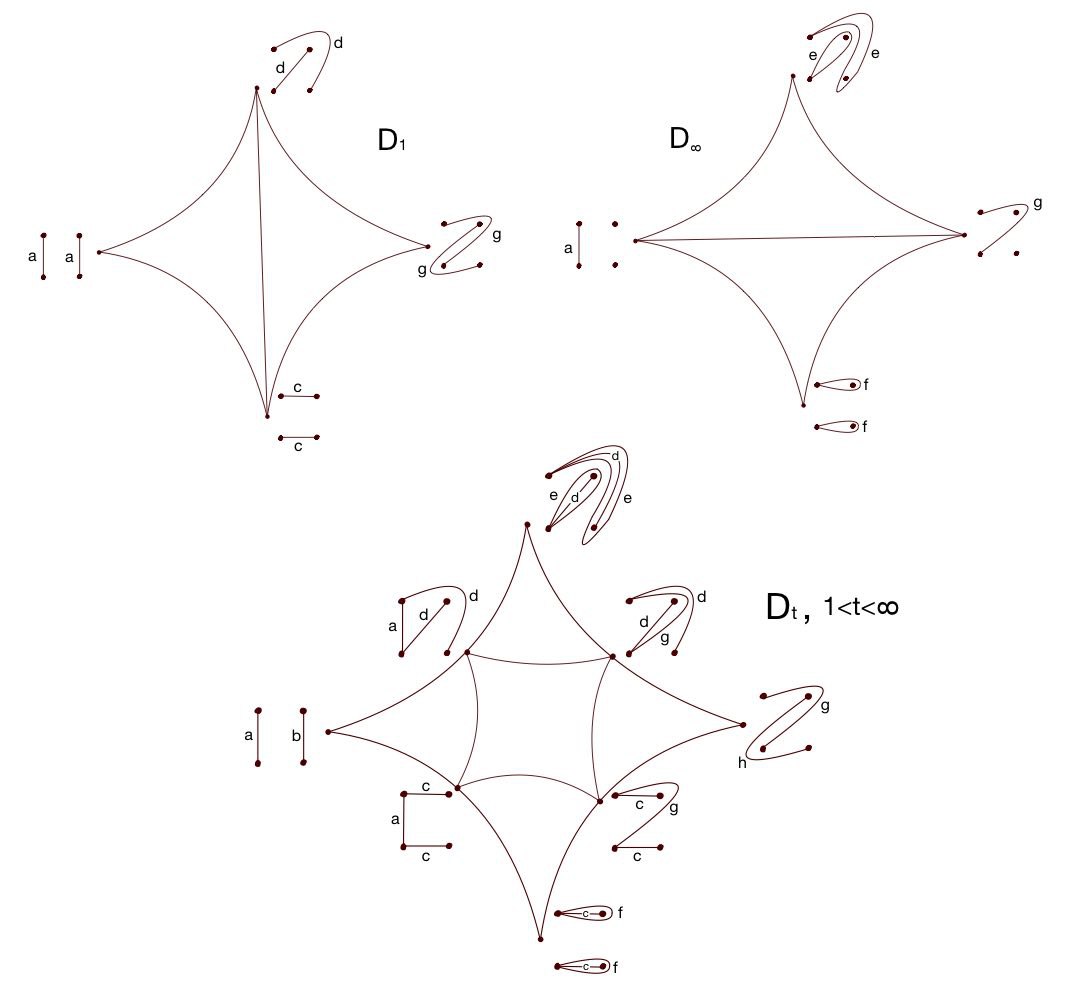}
    \caption{}\label{fig: arc-systems}
\end{figure}

 We now associate to every vertex of $D_\alpha$ a graph in $\Ss$ with vertices in the marked points.
 In Figure \ref{fig: arc-systems} are pictured the arc-systems in the quadrilateral $(1/0,0/1,1/2,1/1)$ of $D_\alpha$. Here, an arc labeled $n$ represents a parallel union of $n$ copies of itself, intersecting just at the endpoints. The group $G\subset \psl (\Z)$ acts on $\Ss=\R^2/\Gamma$ (cf. subsection \ref{subsec: links}) as well. We then extend the assignment of arc-systems in $(1/0,0/1,1/2,1/1)$ to every other quadrilateral in $D_\alpha$ through the action of $G$. Finally, we associate to every edge of $D_\alpha$, the arc-system obtained as the union of the two arc-systems at its endpoints, isotoped to meet just in the marked points.

When $S$ is a surface isotoped as in the hypotesis of Theorem \ref{thm: fh1} (in particular in the case of our interest), then the graph component of $S_t=S\cap\Ss\times \{t\}$ for $t$ close to $1$ is (isotopic to) the arc-system corresponding to the vertex $1/0$ in $D_{a/b}$, for some choice of the labelings $a,b$. Notice that, under the Assumptions \ref{assumptions}, $a$ and $b$ are exactly the integers so that $[S]=a\ell_1+b\ell_2$. The graph component of $S_t$ for $t$ close to $0$ is the arc-system corresponding to the vertex $p/q$ of $D_{a/b}$, appropriately labeled. At any nonsingular level $S_t$, the graph component is either the arc-system associated to a vertex of $D_{a/b}$ or the arc-system associated to an edge of $D_{a/b}$, in both cases appropriately labeled. This last fact is not obvious, and is proved by Floyd and Hatcher by means of the hypoteses on $S$ (see Theorem~\ref{thm: fh1} again).

Let $0<t_1<...<t_k<1$ be the sequence of singular heights for $S$. Observe that $k$ is equal to the number of saddle singularities of the height function, so $k$ is strongly related to the Euler characteristic of $S$. If, for some $i$, $t_i<t,t'<t_{i+1}$, then $S_t$ and $S_{t'}$ are isotopic and their arc-systems are represented by a same vertex or edge in $D_{a/b}$. Indicate with $\lambda_i$ the data given by this arc-system and its labeling. Up to isotopy, we can assume that every $\lambda_i$ is different from $\lambda_{i-1}$, that is we can eliminate every saddle where the arc-system does not change. If $\lambda_i$ is on a vertex, then $\lambda_{i-1}$ and $\lambda_{i+1}$ are on an edge containing that vertex. If $\lambda_i$ and $\lambda_{i+1}$ are on edges, then they are on the same one, but with different labelings. So, the surface $S$ describes a path on the $1$-skeleton of $D_{a/b}$ with endpoints at the vertices $1/0$ and $p/q$. 

\bt\cite[section 7]{fh}\label{thm: fh2} Let $S\subset M$ be a properly embedded oriented surface which is incompressible, meridional-incompressible, and has no peripheral components. Let $\lambda_0,...,\lambda_{k+1}$ be a sequence of arc-systems as above, obtained after isotoping $S$ as in Theorem \ref{thm: fh1}. This sequence satisfies:
\begin{enumerate}
    \item Monotonicity: if $\lambda_{i-1},\lambda_i$ and $\lambda_{i+1}$ lie on the same edge, then, for every fixed arc, its label in $\lambda_{i-1},\lambda_i$ and $\lambda_{i+1}$ is either increasing or decreasing.
    \item Minimality: $\lambda_i\ne \lambda_j$ for $i\ne j$, and two different edges of the same $2$-cell of $D_{a/b}$ cannot be both traversed by the path.
\end{enumerate}
\et

\subsection{The coherently oriented case}
\begin{figure}[ht]
    \centering
    \includegraphics[width=\textwidth]{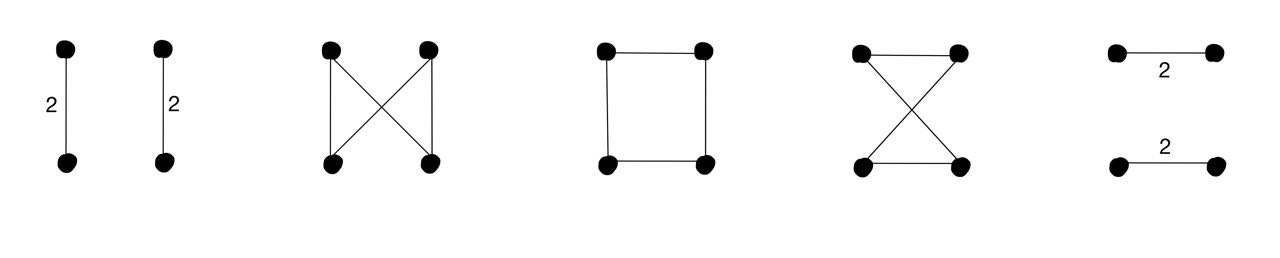}
    \caption{Possible values of $\lambda_0,\lambda_1$ and $\lambda_2$ (first, third and fifth pictures respectively) and their singular levels in the middle, for $S$ representing $2\ell_1+2\ell_2$.}\label{fig: sequences}
\end{figure}

Recall that we are interested in the case when $S$ is a surface satisfying the Assumptions \ref{assumptions}. In that case, the hypothesis of coherent orientation at the boundary imposes strict conditions on the $\lambda$-sequence of $S$. Firstly, if $\lambda_i$ is at a vertex, then the labelings are uniquely determined since every marked point has the same degree in every section $S_t$. Secondly, if $\lambda_i$ is at an edge, then there are finitely many possible labelings for its arc-system, and they are all traversed by the $\lambda$-sequence (for an example, see Figure \ref{fig: sequences}). In particular, the $\lambda$-sequence is recoverable from the path it traverses on $D_{a/b}$. Finally and more importantly, many arc-systems are not realizable anymore, as we will explain.

\begin{figure}[ht]
    \centering
    \includegraphics[width=0.35\textwidth]{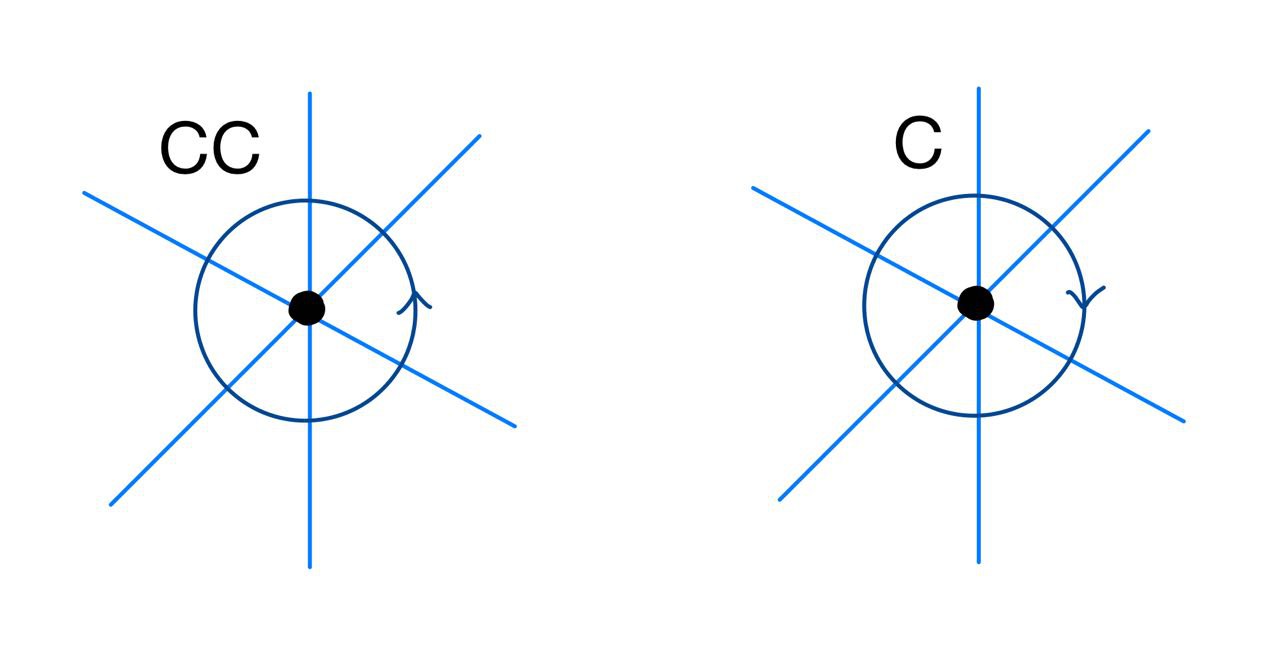}
    \caption{}\label{fig: orientation}
\end{figure}

The surface $S$ is oriented, so it has a normal direction. This orientation induces a normal orientation on every nonsingular section $S_t\subset \Ss\times\{t\}$. The coherent orientation at the boundary allows us assign an orientation to the marked points of $\Ss$. Indeed, we assign a label $C$ to a marked point if and only if every arc pointing at it has the clockwise normal, and a label $CC$ if and only if every arc pointing at it has the counter-clockwise normal, see Figure \ref{fig: orientation}. This is possible because the orientation of $L_{p/q}$ determines the orientation of $S$ close to $\partial S$. Again because of the coherent orientation and of the choice of the orientation so that $a,b\ge 0$, any arc in $S_t$ has a $C$-labeled endpoint and a $CC$-labeled endpoint. We call \emph{admissibly oriented} the arc-systems of $D_{a/b}$ such that, if the marked points are labeled $C$ or $CC$ as in $S_t$ for $t$ close to $1$, each arc of $D_{a/b}$ has a $C$-labeled endpoint and a $CC$-labeled endpoint.

\begin{figure}[ht]
    \centering
    \includegraphics[width=0.8\textwidth]{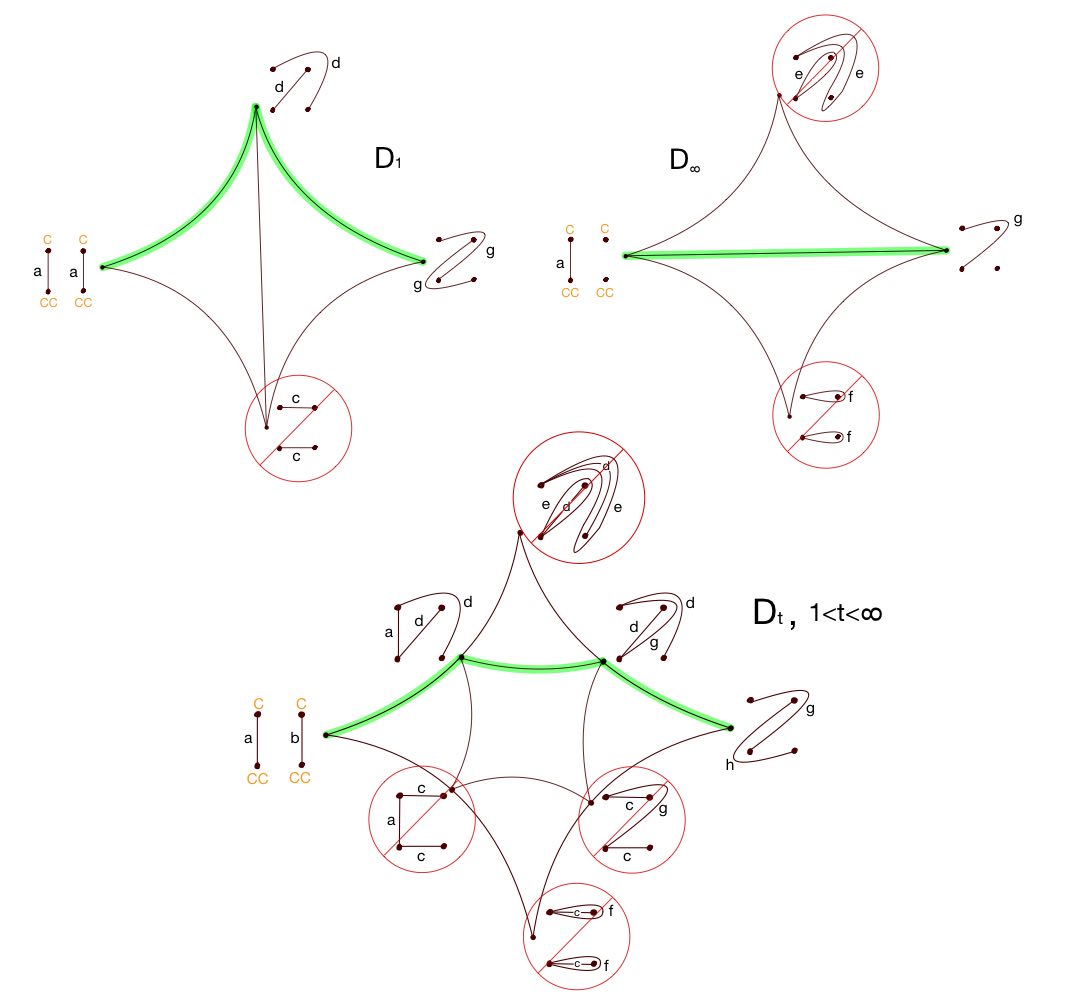}
    \caption{The case when $S$ has normal orientation pointing right, see the labels $C$ and $CC$ at the marked points in $1/0$. Many arc-systems are not admissibly oriented, just the green path can be traversed by a $\lambda$-sequence.}\label{fig: exclusion}
\end{figure}

This condition drastically reduces the set of possible arc-systems in $S_t$. For example, consider Figure \ref{fig: exclusion}, where we suppose that, for $t$ close to $1$, the normal orientation of $S_t$ points right. Then in every quadrilateral of $D_{a/b}$ just half of the vertex arc-systems are admissibly oriented. In general, whatever the normal orientation of $S$ is, on the boundary of every $2$-cell of $D_{a/b}$ there is just one edge and its vertex configurations which are admissibly oriented. 

\bl\label{lemma: tree}
Let $T_{a/b}$ be the subgraph of $D_{a/b}$ composed of the cells of $D_{a/b}$ associated with admissibly oriented arc-systems. Then $T_{a/b}$ is a tree.
\el
\bp Let $\gamma$ be an embedded cycle on the $1$-skeleton of $D_{a/b}$. We show that $\gamma$ cannot be contained in $T_{a/b}$. First suppose that $a/b=0,1$ or $\infty$. Let $e$ be an edge of $\gamma$ and let $A_1$ be a closed triangle of $D_{a/b}$ containing it. Inductively, define $A_i$ from $A_{i-1}$ as the subcomplex of $D_{a/b}$ obtained by adding to $A_{i-1}$ the closed triangles adjacent to it. There is a natural $n$ such that $\gamma$ is not contained in $A_n$ but is contained in $A_{n+1}$. Since $T_{a/b}$ contains just an edge for every $2$-cell of $D_{a/b}$, it is easily seen that $\gamma$ cannot sit inside $T_{a/b}$.

For general values of $a/b$, the same line of reasoning works. In this case, though, we choose $A_1$ to be a quadrilateral tile, i.e. the union of a quadrilateral $2$-cell and its four adjacent triangular $2$-cells, and then build $A_{i+1}$ from $A_i$ by adding all the neighboring quadrilateral tiles instead of just the triangles.
\ep

\bigskip

\bprop\label{prop: admissibly oriented} Let $S\subset M$ be a surface representing the class $a\ell_1+b\ell_2$ ($a\ge b\ge 0$) and satisfying the Assumptions \ref{assumptions}, namely $S$ being properly embedded, oriented, incompressible, meridionally-incompressible, coherently oriented at the boundary (apart for meridians) and with no closed or peripheral components.

Isotope $S$ as in Theorem \ref{thm: fh1} and consider its $\lambda$-sequence $\lambda_0,...,\lambda_{k+1}$. Then:
\begin{enumerate}
    \item The $\lambda$-sequence induces the unique path in $T_{a/b}$ with endpoints $1/0$ and $p/q$;
    \item The surface $S$ is norm-minimizing and $$x(a\ell_1+b\ell_2)=k-a-b.$$
\end{enumerate}
\eprop
\bp The previous lemma and Theorem \ref{thm: fh2} easily imply $(1)$.

To show $(2)$, observe that we can use the $\lambda$-sequence and the singular sections $S_{t_i}$ to construct a handle-decomposition of $S$. Indeed, $S\cap \Ss\times [1-\epsilon,1]$ is the union of $a+b$ $2$-cells ($0$-handles). Then, the surface $S\cap \Ss\times [t_i-\epsilon,1]$ can be obtained from $S\cap \Ss\times [t_{i+1}-\epsilon,1]$ by just gluing one $1$-handle along its boundary. Attaching a $1$-handle decreases the Euler characteristic by $1$, so $\chi(S)=a+b-k$.

Finally, given a norm-minimizing representative of $a\ell_1+b\ell_2$ satisfying the conditions listed at the beginning of subsection \ref{subsec: good representative} and isotoped as in Theorem \ref{thm: fh1}, its $\lambda$-sequence induces the same unique path on the tree $T_{a/b}$ as the one of $S$. In particular, $S$ has the same Euler characteristic of a norm-minimizing surface, thus $S$ is norm-minimizing as well.
\ep

\subsection{The main theorem}

\bt\label{thm: Thurston ball}
Let $L_{p/q}$ be an oriented $2$-bridge link with components $\ell_1$ and $\ell_2$, and let $M=S^3-\overset{\circ}N(L_{p/q})$ be the link exterior.
   The Thurston unit ball in $H_2(M,\partial M;\R)$ is the convex polygon spanned by the points $\pm\frac{\ell_1}{x(\ell_1)}$, $\pm\frac{\ell_2}{x(\ell_2)}$, $\pm\frac{\ell_1+\ell_2}{x(\ell_1+\ell_2)}$, and $\pm\frac{\ell_1-\ell_2}{x(\ell_1-\ell_2)}$.
\et
\bp
Since there is a self-homeomorphism of $S^3$ exchanging the components of $L_{p/q}$, and thanks to Lemma \ref{lemma: additivity}, it is enough to show that $x(2\ell_1\pm\ell_2)=x(\ell_1)+x(\ell_1\pm\ell_2)$. The key observation is that, for $a/b=1/1,2/1$ and $1/0$, the only arc-systems of $D_{a/b}$ are on vertices. Indeed, starting from an arc-system relative to a vertex, every permitted saddle transforms the system into an arc-system at another vertex, because of the small number of arcs. Let $\gamma_2$ be the embedded path in $T_{2/1}$ with endpoints $1/0$ and $p/q$. Thanks to Proposition \ref{prop: admissibly oriented}: $$x(2\ell_1+\ell_2)=\lg(\gamma_2)-3$$ where $\lg(\gamma_2)$ is the number of edges of $\gamma_2$.

Every admissibly oriented arc-system in $D_{2/1}$ is the disjoint union of an admissibly oriented arc-system in $D_{1/1}$ and one in $D_{1/0}$. Call $\lambda_0^{(2)},...,\lambda_k^{(2)}$ the sequence of admissible points defining $\gamma_2$, and consider for every $i$ the configurations $\lambda_i^{(1)}$ on $ D_{1/1}$ and $\lambda_i^{(\infty)}$ on $ D_{1/0}$ whose arc-systems' union gives the arc-system of $\lambda_i^{(2)}\in D_{2/1}$. Observe that $\lambda_i^{(1)}=\lambda_{i+1}^{(1)}$ if and only if $\lambda_i^{(\infty)}$ and $\lambda_{i+1}^{(\infty)}$ are on different endpoints of an edge of $D_{1/0}$. Also  $\lambda_i^{(\infty)}=\lambda_{i+1}^{(\infty)}$ if and only if $\lambda_i^{(1)}$ and $\lambda_{i+1}^{(1)}$ are on different endpoints of an edge of $D_{1/1}$. In particular, $$\lg(\gamma_2)=\lg(\gamma_1)+\lg(\gamma_{\infty}),$$ where $\gamma_1$ is the path in $D_{1/1}$ passing through $\lambda_0^{(1)},...,\lambda_k^{(1)}$ and $\gamma_{\infty}$ is the path in $D_{1/0}$ passing through $\lambda_0^{(\infty)},...,\lambda_k^{(\infty)}$. 

Finally, we get \begin{align*}
x(2\ell_1+\ell_2)&=\lg(\gamma_2)-3=\lg(\gamma_1)-2+\lg(\gamma_{\infty})-1= \\ &= x(\ell_1+\ell_2)+x(\ell_1).    
\end{align*}
An analogous reasoning shows that $x(2\ell_1-\ell_2)=x(\ell_1)+x(\ell_1-\ell_2)$.
\ep

\br\label{rem: cut-and-paste} The relation $\lg(\gamma_2)=\lg(\gamma_1)+\lg(\gamma_{\infty})$ has a geometric counterpart. Consider a surface $S_2$ realized by the path $\gamma_2$. We claim that $S_2$ is the oriented cut-and-paste of two surfaces $S_1$ and $S_{\infty}$, respectively inducing the paths $\gamma_1$ and $\gamma_{\infty}$, and so representing the classes $\ell_1+\ell_2$ and $\ell_1$. Close to a nonsingular level $t$, the surface $S_2$ is just the product $(S_2)_t\times [-\epsilon,\epsilon]$. At the same time, $(S_2)_t$ is uniquely decomposable as the disjoint union of an arc-system in $D_{1/1}$ and one in $D_{1/0}$. So, close to height $t$, we define $S_1$ as the product of the arc-system in $D_{1/1}$ times $[-\epsilon,\epsilon]$, and analogously for $S_{\infty}$.

Let us consider what $S_2$ looks like near a singular level $t$. Here, the behavior is determined by an edge $e$ of $\gamma_2$. Each arc-system at the endpoints of $e$ can be uniquely decomposed as the disjoint union of an arc-system in $D_{1/1}$ and one in $D_{1/0}$. The edge $e$ is collapsed to a point in $\gamma_1$ if and only if the two arc-systems at the endpoints of $e$ have the same summand coming from $D_{1/1}$. Otherwise, $e$ is collapsed to a point in $\gamma_{\infty}$ and the two arc-systems at the endpoints of $e$ have the same summand coming from $D_{1/0}$. We consider these two possibilities separately.

\begin{figure}
    \centering
    \includegraphics[width=0.8\textwidth]{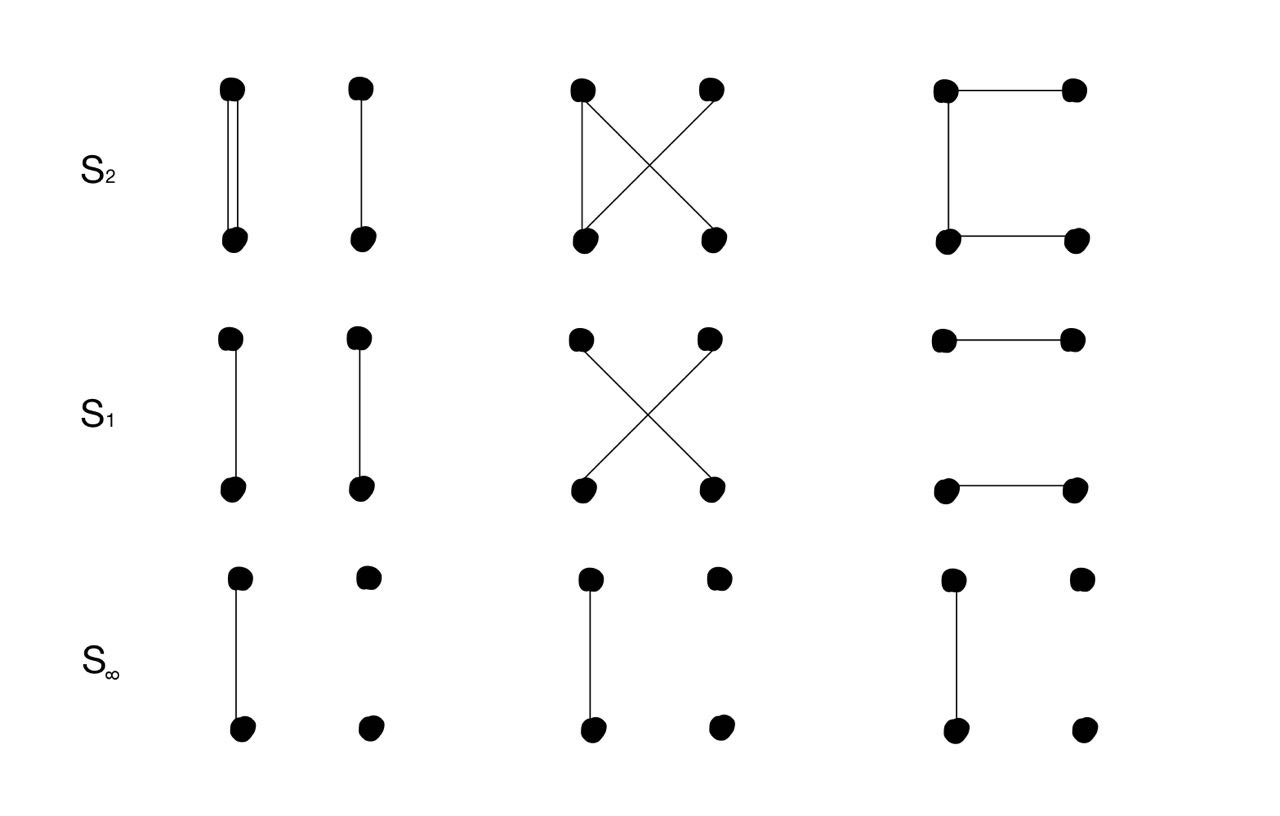}
    \caption{From left to right: the sections of the surfaces $\phi(S_2),\phi(S_1)$ and $\phi(S_\infty)$ respectively at heights $t-\epsilon, t$ and $t+\epsilon$.}\label{fig: sections1}
\end{figure}

If $e$ is collapsed to a point in $\gamma_\infty$, then there is a self-diffeomorphism $\phi$ of $S^2_*\times[t-\epsilon,t+\epsilon]$ sending $(S_2)_{t-\epsilon}$, $(S_2)_{t}$ and $(S_2)_{t+\epsilon}$ to the sections in Figure \ref{fig: sections1} top row, possibly in reversed order. In this figure, the sections definining $\phi(S_1)$ and $\phi(S_{\infty})$ are illustrated. 
Also in this case, $S_2$ is just the disjoint union of $S_1$ and $S_{\infty}$.

\begin{figure}
    \centering
    \includegraphics[width=0.8\textwidth]{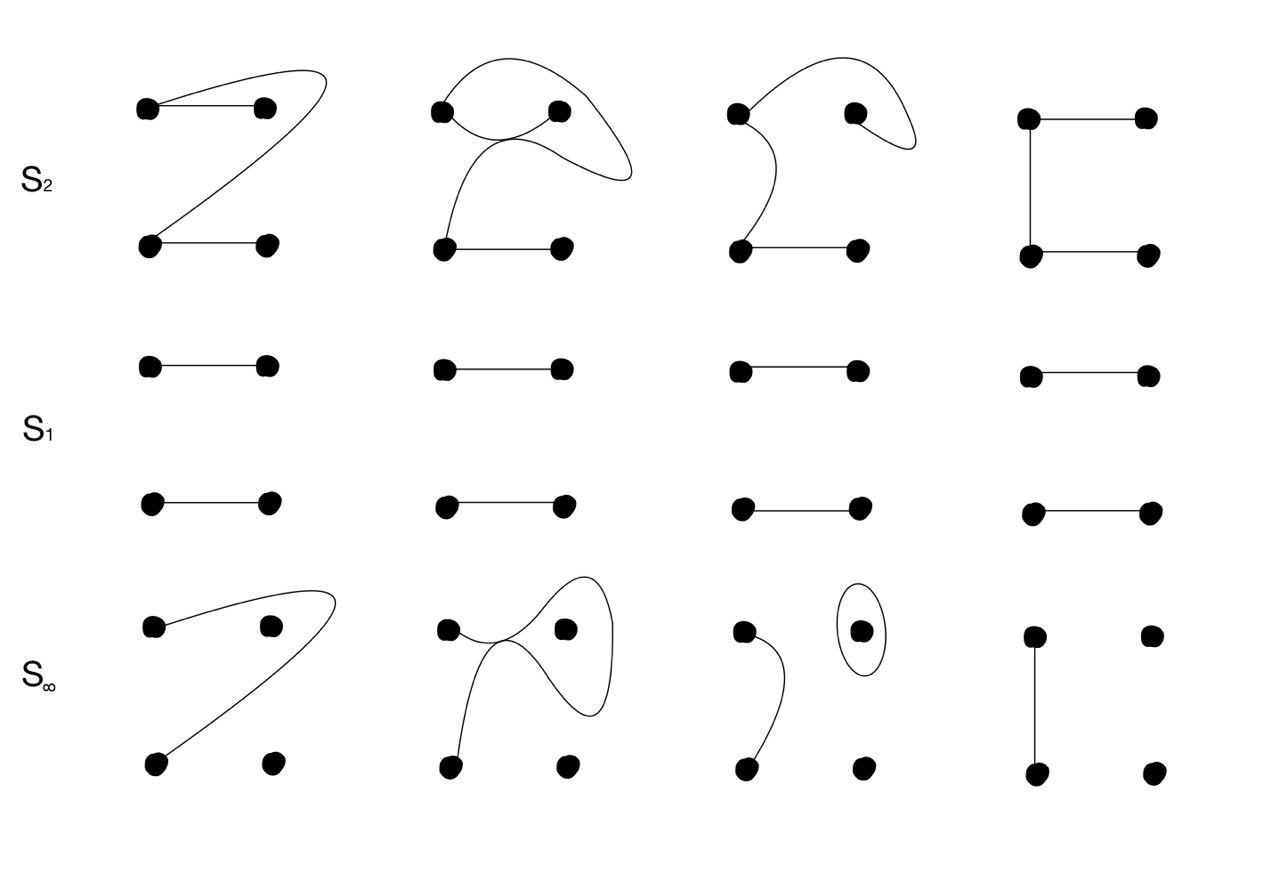}
    \caption{From left to right: the sections of the surfaces $\phi(S_2),\phi(S_1)$ and $\phi(S_\infty)$ respectively at heights $t-\epsilon, t, t+\epsilon/2$ and $t+\epsilon$.}\label{fig: sections2}
\end{figure}

\begin{figure}
    \centering
    \includegraphics[width=.95\textwidth]{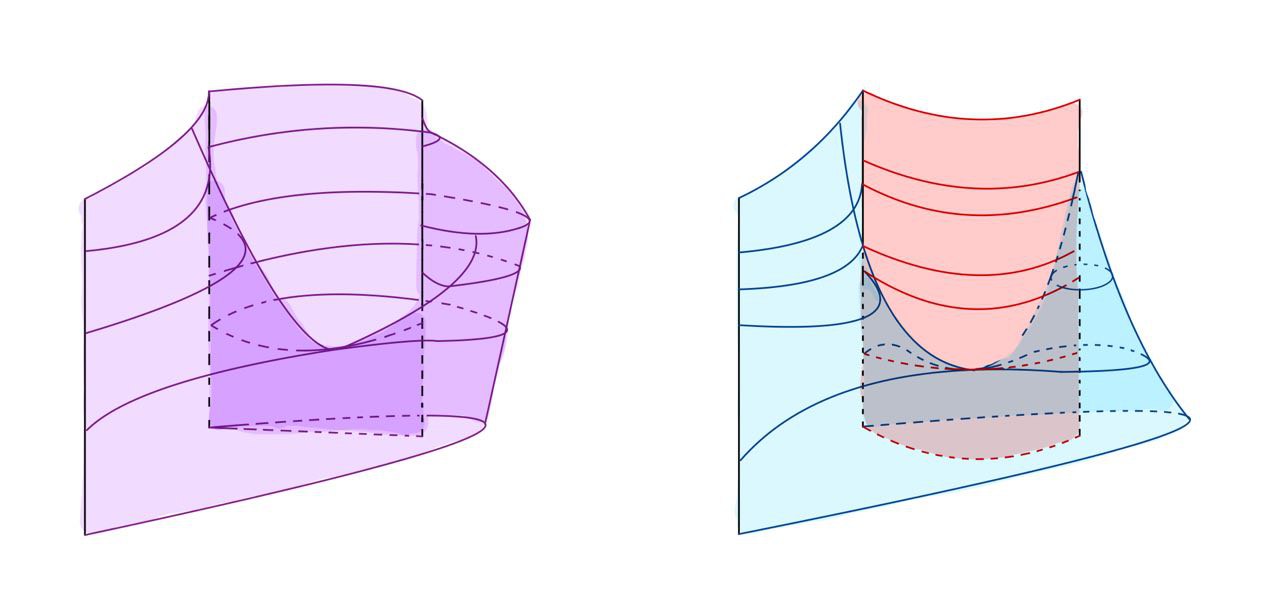}
    \caption{The surface $\phi(S_2)$ near a saddle (left) is isotopic to the oriented cut-and-paste of $\phi(S_1)$ (right, red) and $\phi(S_{\infty})$ (right, light blue). In both figures, a product rectangle should come from the bottom horizontal arc of Figure \ref{fig: sections2} first and second rows. We are deleting this rectangle so as not to complicate the picture anymore.}\label{fig: cutandpaste}
\end{figure}

If $e$ is collapsed to a point in $\gamma_1$, then there is a self-diffeomorphism $\phi$ of $S^2_*\times[t-\epsilon,t+\epsilon]$ sending $(S_2)_{t-\epsilon}$, $(S_2)_{t}$, $(S_2)_{t+\epsilon/2}$ and $(S_2)_{t+\epsilon}$ to the sections in Figure \ref{fig: sections2} top row, possibly in reversed order. In this figure, the sections definining $\phi(S_1)$ and $\phi(S_{\infty})$ are illustrated. 
In this case, we can easily visualize that $S_2$ is the oriented cut-and-paste of $S_1$ and $S_{\infty}$ (see Figure \ref{fig: cutandpaste}).

\er

\newpage
\section{Norm-minimizing surfaces}\label{sec: norm-minimizing}
The goal of this section is to exhibit norm-minimizing representatives for integral homology classes in $2$-bridge link complements. This will give us an algorithmic way to compute the Thurston norm of any homology class.

\subsection{Base-type diagrams} 
Base-type $2$-bridge links can be characterized as the $2$-bridge links whose Thurston ball is a quadrilateral with vertices along the bisectors or as the ones whose link complements fibers in a specific way over $S^1$ (see Corollary \ref{cor: base-type}). Every $2$-bridge link is a tangle sum of base-type $2$-bridge links, and this tangle sum happens to behave well with respect to the Thurston norm. So, base-type $2$-bridge links are the somewhat natural building blocks for the computation of the Thurston norm of $2$-bridge link complements. 

\begin{definition}
    Let $L=T(a_1,a_2,...,a_k)$ be an alternating diagram for a $2$-bridge link. We say that $L$ is \emph{base-type}, if either $k=1$ and $a_1$ is even, or if $k>1$, $a_1$ and $a_k$ are odd, whilst $a_i$ is even for every $1<i<k$.
\end{definition}

\br

\begin{figure}[ht]
    \centering
    \includegraphics[width=.9\textwidth]{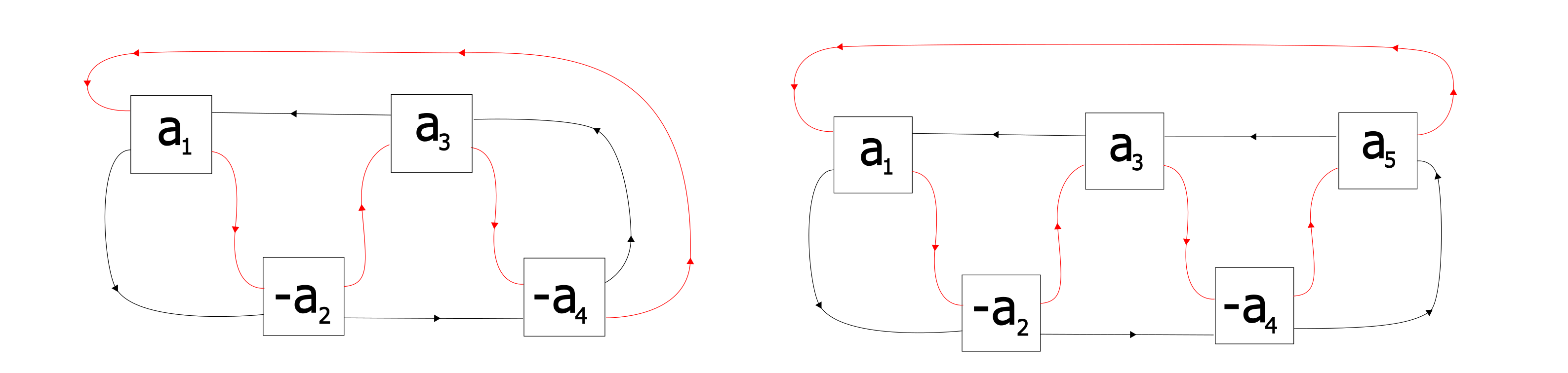}
    \caption{Two examples of base-type diagrams with their standard orientation.}\label{fig: base type}
\end{figure}

    \begin{itemize}
        \item Base-type diagrams are the alternating rational diagrams of $2$-bridge links such that no component of the link crosses itself. See Figure~\ref{fig: base type}. Notice that it is particularly easy to see how base-type diagrams are standardly oriented.
        \item If a $2$-bridge link has a base-type diagram, then every alternating rational diagram representing it is base-type. So, we will refer to these links as \emph{base-type} $2$-bridge links.
    \end{itemize}
\er

\begin{lemma}
    Let $L$ be an oriented nonsplit link with components $\ell_1$ and $\ell_2$, and let $M=M_L$ be the link exterior in $S^3$. Suppose that there is a minimal-genus surface $D$ for $\ell_1$ in $S^3$ such that $\ell_2$ intersects $D$ transversely and always with the same sign. 
    Then the surface $D\cap M$ is norm-minimizing in $H_2(M,\partial M;\R)$.
\end{lemma}

\bp Identify $H_2(M,\partial M;\Z)\cong H_1(\ell;\Z)=(\Z \ell_1) \oplus (\Z\ell_2)$. The surface $D\cap M$ is obtained from $D$ by removing $|D\cap\ell_2|=|i(D,\ell_2)|=|\lk(\ell_1,\ell_2)|$ discs, and it represents the class $\ell_1$ in $H_2(M,\partial M;\Z)$ (here $i(X,Y)$ represents the algebraic intersection between $X$ and $Y$). Another surface $S\subset M$ representing $\ell_1$ must intersect $\partial N(\ell_2)$ in meridional loops. Thus $S=\overline S\cap M$, where $\overline S\subset S^3$ is a Seifert surface for $\ell_1$. Observe that since $\overline S$ is a Seifert surface, $\chi(\overline S)\ge 0$ if and only if $\chi(\overline S)=1$. In particular, $\chi_-(D)\le\chi_-(\overline S)$ implies $\chi(D)\ge\chi(\overline S)$. Now we have $$\chi(S)=\chi(\overline S)-|\overline S\cap \ell_2|\le \chi(\overline S)-|i(S,\ell_2)|=\chi(\overline S)-|\lk(\ell_1,\ell_2)|\le \chi(D)-|D\cap\ell_2|=\chi(D\cap M).$$
\ep 

\bc\label{cor: base type 10}
    Let $L=T(a_1,a_2,...,a_k)$ be a base-type diagram for a $2$-bridge link with components $\ell_1$ and $\ell_2$. By identifying $H_2(M_L,\partial M_L;\Z)=(\Z \ell_1) \oplus (\Z\ell_2)$ we have $$x(\ell_1)=x(\ell_2)=\frac 12\sum_{i=1}^k |a_i|-1.$$
\ec

\bp The projections of $\ell_1$ and $\ell_2$ in $T(a_1,...,a_k)$ bound a disc which always intersects the other component with the same sign. So, the previous lemma applies.
\ep

\begin{figure}[h]
    \centering
    \includegraphics[width=0.6\textwidth]{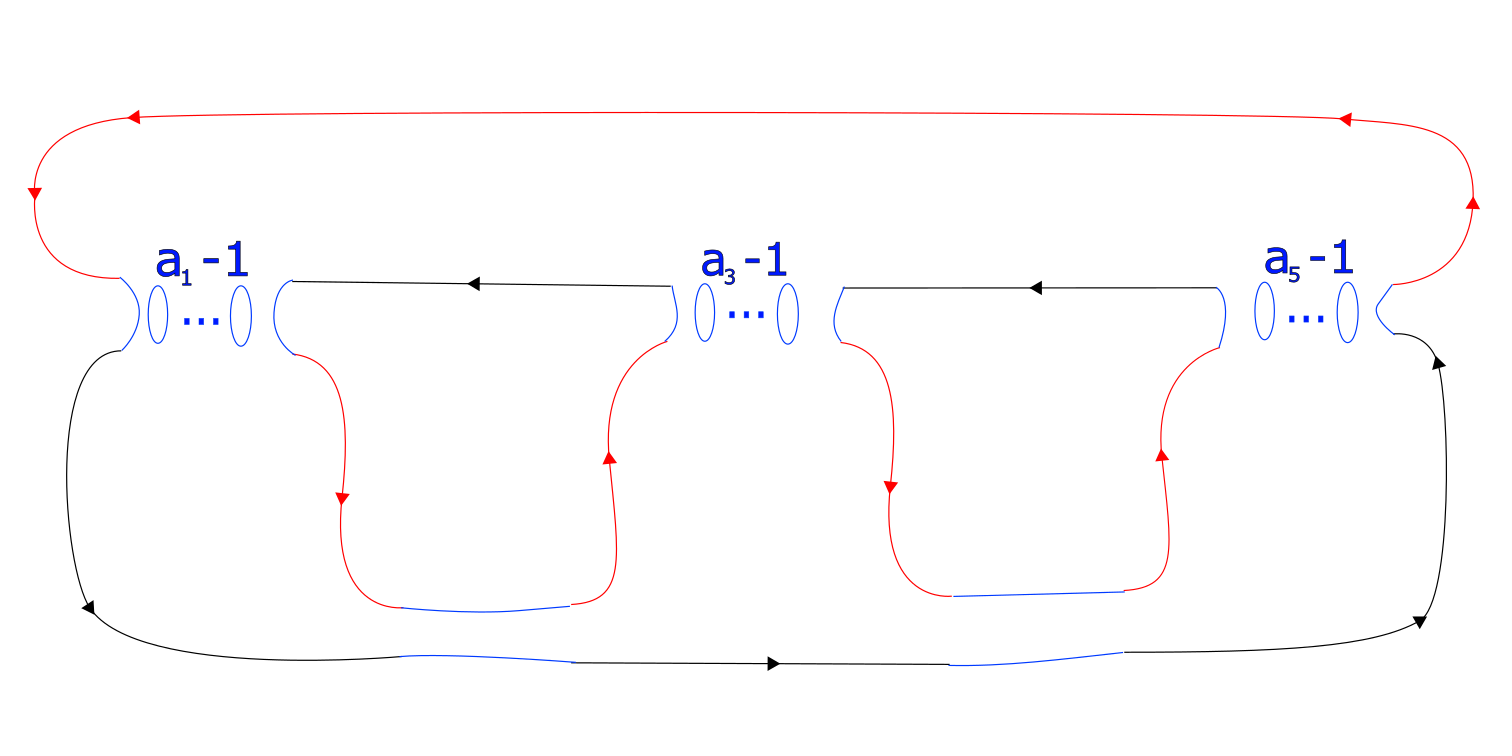}
    \caption{Seifert circles for a base-type diagram of odd length.}\label{fig: seifert}
\end{figure}

\bprop\label{prop: base type 11}
    Let $L=T(a_1,a_2,...,a_k)$ be a standardly oriented base-type diagram for a $2$-bridge link with components $\ell_1$ and $\ell_2$. By identifying $H_2(M_L,\partial M_L;\Z)=(\Z \ell_1) \oplus (\Z\ell_2)$ we have $$x(\ell_1+\ell_2)=\begin{cases}
    \sum_{i=1}^{\frac{k-1}2} |a_{2i}|\; \; \text{if $k$ is odd} \\ \\
    \sum_{i=1}^{\frac{k}2} |a_{2i}|-1 \; \; \text{if $k$ is even} 
    \end{cases}$$
    $$x(\ell_1-\ell_2)=\begin{cases}
    \sum_{i=0}^{\frac{k-1}2} |a_{2i+1}|-2\; \; \text{if $k$ is odd} \\ \\
    \sum_{i=0}^{\frac{k}2-1} |a_{2i+1}|-1 \; \; \text{if $k$ is even} .
    \end{cases}$$
\eprop
\bp Since $L$ is an alternating diagram, applying Seifert algorithm to $L$ gives a norm minimizing surface for $\ell_1+\ell_2$ (cf. \cite{murasugi, crowell}). For odd values of $k$, the algorithm produces $|a_i|-1$ circles in every box with odd $i$, one circle over every box with even $i$, all nested in another circle given by the ``external perimeter" of the link diagram (see Figure~\ref{fig: seifert}). Thus, the obtained Seifert surface has Euler characteristic $$\left(\sum_{i=0}^{\frac{k-1}2} (|a_{2i+1}|-1)+\frac{k-1}2+1\right)-\sum_{i=1}^{k} |a_{i}|=-\sum_{i=1}^{\frac{k-1}2} |a_{2i}|.$$ 
The other computations are analogous.
\ep

\bc\label{cor: base type}
    Let $L$ be a base-type oriented $2$-bridge link and let $M_L$ be its exterior in $S^3$. The Thurston norm $x$ on $H_2(M_L,\partial M_L;\R)$ is the polygon spanned by the points $\pm\frac{\ell_1+\ell_2}{x(\ell_1+\ell_2)}$ and $\pm\frac{\ell_1-\ell_2}{x(\ell_1-\ell_2)}$.
\ec
\bp By the previous computations, $x(\ell_1)=x(\ell_2)=\frac 12 \left(x(\ell_1+\ell_2)+x(\ell_1-\ell_2)\right)$. Then, Lemma \ref{lemma: additivity} concludes.\ep

\subsection{The general case}
Let $L=T(a_1,...,a_k)$ be an alternating rational diagram. Endow $L=\ell_1\sqcup \ell_2$ with the standard orientation, and identify $H_2(M_L,\partial M_L;\Z)=(\Z \ell_1) \oplus (\Z\ell_2)$. Now, we define a standard representative for each class $a\ell_1+b\ell_2$. 

\begin{figure}[ht]
    \centering
    \includegraphics[width=0.6\textwidth]{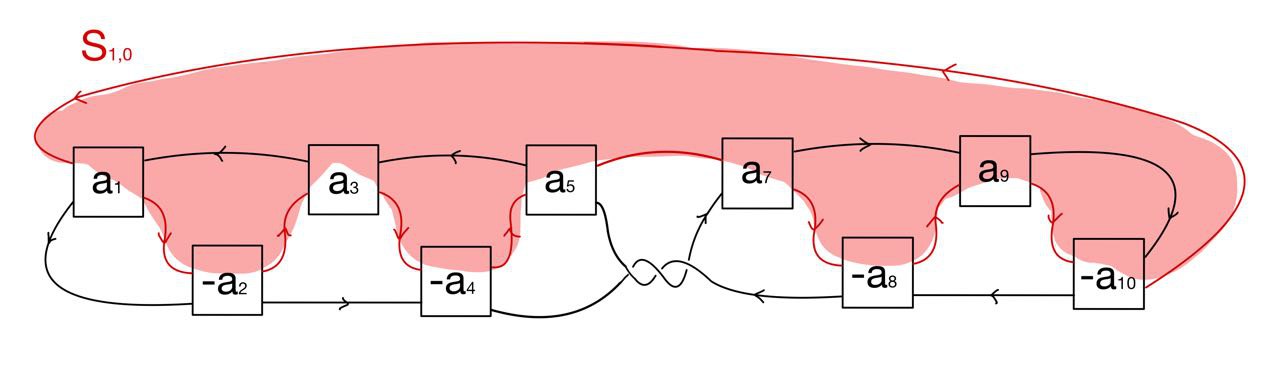}
    \includegraphics[width=0.6\textwidth]{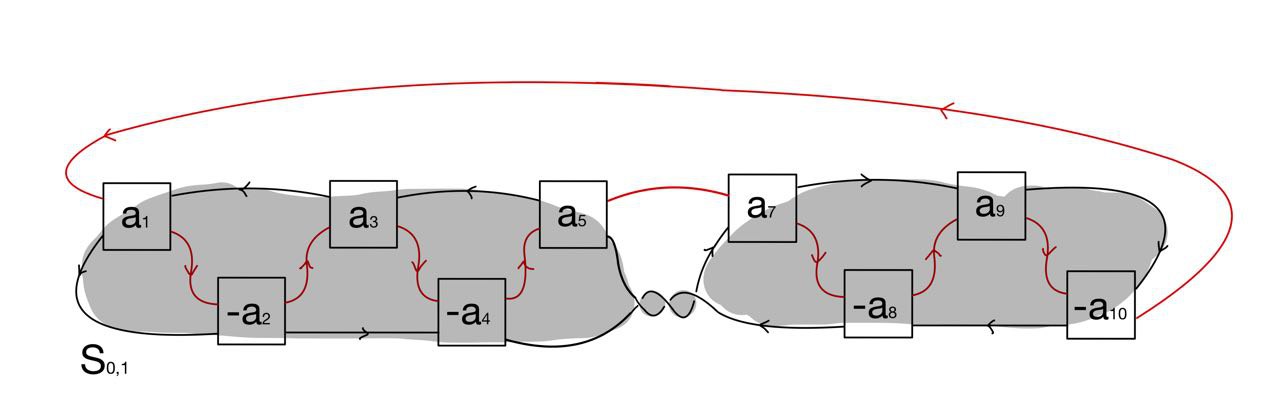}
    \caption{}\label{fig: standard representatives}
\end{figure}

Let $S_{1,0}(L)$ be the surface $D\cap M_L$, where $D\subset S^3$ is the oriented disc that projects to the bounded disc with boundary $\ell_1$ in $L\subset \R^2$. Similarly, let $S_{0,1}(L)$ be the surface $D'\cap M_L$, where $D'\subset S^3$ is the oriented disc that projects to the bounded twisted disc with boundary $\ell_2$ in $L\subset \R^2$. See Figure~\ref{fig: standard representatives}.

Let $S_{1,\pm 1}(L)$ be the surfaces obtained from Seifert algorithm applied to the diagram $L$ endowed with the relative orientation, where we choose the bounded disc for each Seifert circle.

If $a\ge b\ge 0$, let $S_{a,b}(L)$ be the oriented cut-and-paste sum of $a-b$ copies of $S_{1,0}(L)$ and of $b$ copies of $S_{1,1}(L)$. Let $S_{a,-b}(L)$ be the oriented cut-and-paste sum of $a-b$ copies of $S_{1,0}(L)$ and of $b$ copies of $S_{1,-1}(L)$.

If $b\ge a\ge 0$, define $S_{a,\pm b}(L)$ analogously by means of $S_{0,1}$ instead of $S_{1,0}$.

Finally, if $a\le 0$, define $S_{a,b}(L)$ just as $S_{-a,-b}(L)$ with the opposite orientation.

Clearly, $S_{a,b}(L)$ represents the class $a\ell_1+b\ell_2$ in $H_2(M_L,\partial M_L;\Z)$.

\begin{figure}[ht]
    \centering
    \includegraphics[width=0.8\textwidth]{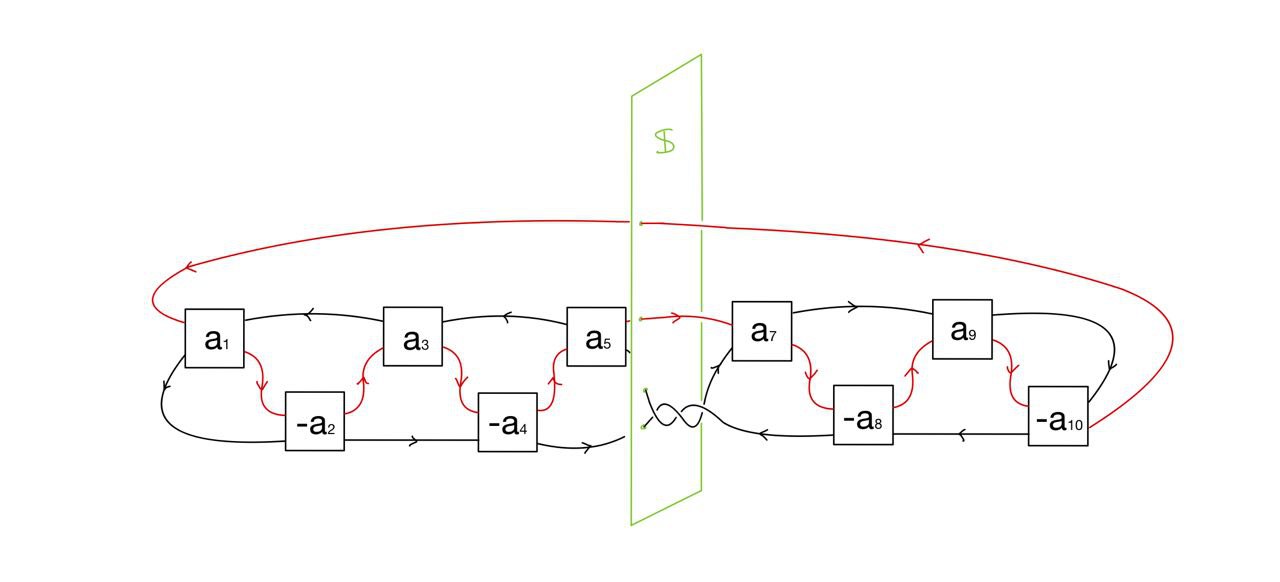}
    \caption{}\label{fig: sphere}
\end{figure}

\br\label{rem: tangle} If $L$ is not base-type, then there is an $i$ such that the component $\ell_2$ crosses itself in the $i$-th box. Let $\mathbb S\subset S^3$ be a sphere intersecting $L$ transversely $4$ times and in the $i$-th box, as in Figure~\ref{fig: sphere} (where $i=6$).

The surface $S_{a,b}(L)$ intersects $\mathbb S$ transversely in $|a|$ parallel arcs connecting two points of $\ell_1$ and $|b|$ parallel arcs connecting two points of $\ell_2$.

\begin{figure}[ht]
    \centering
    \includegraphics[width=0.8\textwidth]{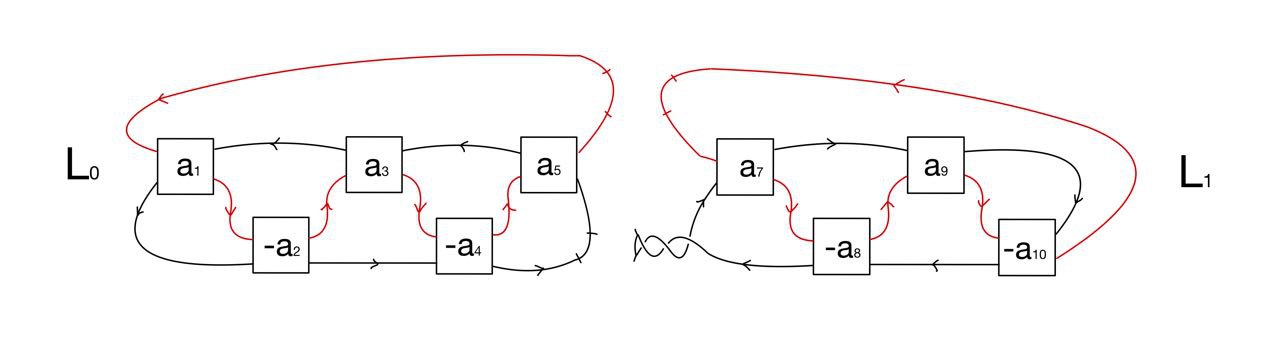}
    \caption{}\label{fig: tangle sum}
\end{figure}

Cut the link $L$ along $\mathbb S$ and glue on each side a copy of one of both types of arcs coming from $S_{a,b}(L)\cap\mathbb S$, so to obtain two $2$-bridge diagrams $L_0$ and $L_1$. See Figure~\ref{fig: tangle sum}. We refer to $L$ as a \emph{tangle sum} of $L_0$ and $L_1$ through $\mathbb S$.

The diagrams $L_0$ and $L_1$ are equivalent to the diagrams $T(a_1,...,a_{i-1})$ and $T(a_{i+1},...,a_k)$ respectively. Note that the induced orientation on $L_1$ is standard if and only if $a_i$ is even.

Now, $S_{a,b}(L)$ can be obtained from $S_{a,b}(L_0)$ and $S_{a,\pm b}(L_1)$ (the sign depending on the parity of $a_i$) by attaching $|a|+|b|$ $1$-handles along their boundaries. In particular \begin{align}\chi(S_{a,b}(L))=\chi(S_{a,b}(L_0))+\chi(S_{a,\pm b}(L_1))-|a|-|b|.\label{eq: tangel sum}\end{align}
\er 

\bt\label{thm: standard representatives} Let $L=T(a_1,...,a_k)$ be an alternating diagram of a $2$-bridge link. For every $a,b\in \Z$, the surface $S_{a,b}(L)$ is a norm-minimizing representative for the class $a\ell_1+b\ell_2$.
\et
\bp Thanks to Theorem~\ref{thm: Thurston ball}, we just need to show that the surfaces $S_{a,b}(L)$ are norm-minimizing for $(a,b)=(1,\pm 1),(1,0)$ and $(0,1)$. The cases $(a,b)=(1,\pm 1)$ are easily proved: Seifert algorithm gives minimal-genus surfaces when applied to alternating diagrams (cf. \cite{murasugi, crowell}), and it is straightforward to verify that the tangle sum of Remark~\ref{rem: tangle} respects Seifert algorithms of the summands. Though, an argument similar to the following one for $(1,0)$ works as well for $(1,\pm 1)$. By the symmetry of $L$, we just need to show that $S_{1,0}$ is norm-minimizing.

We proceed by induction on the number of self-crossing boxes in $T(a_1,...,a_k)$, Corollary \ref{cor: base type 10} giving the base of the induction. Suppose that $a_i$ represents a self-crossing box in $T(a_1,...,a_k)$, dividing it into alternating $2$-bridge diagrams $L_0=T(a_1,...,a_{i-1})$ and $L_1=T(a_{i+1},...,a_k)$. For $j=1,...,k$, define \begin{align*}
    p_j/q_j=[a_1,...,a_j], \\ \bar p_j/\bar q_j=[a_j,...,a_k].
\end{align*} 
Hence, $L_0=L_{p_{i-1}/q_{i-1}}$ and $L_1=L_{\bar p_{i+1}/\bar q_{i+1}}$.
Suppose that the statement of the theorem holds for $L_0$ and $L_1$. Recall that by Lemma \ref{lemma: tree} there is only one embedded path $\gamma$ in $T_{1/0}$, going from $1/0$ to $p/q$. We want to show that $\chi(S_{1,0}(L))=1-\lg(\gamma)$ and we know an analogous statement holds for $L_0$ and $L_1$. Let $\gamma_0$ be the unique embedded path in $T_{1/0}$ going from $1/0$ to $p_{i-1}/q_{i-1}$, and let $\gamma_1$ be the unique embedded path in $T_{1/0}$ going from $1/0$ to $\bar p_{i+1}/\bar q_{i+1}$. Then, for $i=0,1$, $$\chi(S_{1,0}(L_i))=1-\lg(\gamma_i).$$ Thus, thanks to identity (\ref{eq: tangel sum}), we just need to show that $$\lg(\gamma)=\lg(\gamma_0)+\lg(\gamma_1).$$
The group $\SL(\Z)$ acts on $\Z^2=\Q\cup\{\infty\}$ by left-multiplication, and this induces a cellular action of $\psl(\Z)$ on $D_{1/0}$. 

By induction on $j$, one can easily show that \begin{equation}\label{eq: recurrence}
   \begin{cases} p_{j+1}=a_{j+1}p_j+p_{j-1} \\
    q_{j+1}=a_{j+1}q_j+q_{j-1}.
    \end{cases}
\end{equation}
Then, the matrix $T=\begin{pmatrix}
    p_{i-1} & p_{i} \\ q_{i-1} & q_{i}
\end{pmatrix}$ has determinant $1$, and thus belongs to $\SL(\Z)$. Again by induction on $n-i$, the fraction $\bar p_{i+1}/\bar q_{i+1}$ is proved to be sent to $p/q$ by $T$. The admissibly oriented arc-systems in $D_{1/0}$ are the ones associated to reduced fractions with even denominators. Since $T$ sends reduced fractions with even denominators to other reduced fractions with even denominators, $T$ sends admissibly oriented arc-systems to admissibly oriented arc-systems.

Consequently, $T(\gamma_1)$ is an embedded path in $T_{1/0}$ going from $p_{i-1}/q_{i-1}$ to $p/q$, and the concatenation of $\gamma_0$ and $T(\gamma_1)$ gives a path in $T_{1/0}$, going from $1/0$ to $p/q$. That concatenation is also embedded: since all the $a_j$ have the same sign, one can see that the edge between $p_{i-1}/q_{i-1}$ and $p_i/q_i$  (which is not admissibly oriented) separates $\gamma_0$ and $T(\gamma_1)$ in $D_{1/0}$. We conclude that $$\lg(\gamma)=\lg(\gamma_0)+\lg(T(\gamma_1))=\lg(\gamma_0)+\lg(\gamma_1).$$
\ep

\subsection{Thurston norm in $2$-bridge link complements}\label{subsec: simpler}
Thanks to the results of the previous sections, it is now easy to compute the Thurston norm of any second-homology class for a $2$-bridge link complement. Indeed, the norm of any integer class is computed as a linear combination of the norms of the vertices, and the latter are computed by an iterated use of equality (\ref{eq: tangel sum}) after decomposing the alternating rational diagram in base-type pieces.

Now, whenever three of the vertices of Theorem \ref{thm: Thurston ball} happen to be aligned, a fewer number of vertices spans the Thurston ball, and we would expect such a lower complexity of the norm to come from a lower complexity of the link itself. The next results show that this is indeed the case, but we need to introduce some notation first.

\bd Given an oriented $2$-component link $L=\ell_1\cup\ell_2$, consider the plane $H_2(M_L,\partial M_L;\R)$ with basis $\ell_1,\ell_2$ and the Thurston norm $x_L$. We denote by $\rays(L)$ the set of slopes $\alpha\in \Q\cup\{\infty\}$ such that there is a vertex of $B_{x_L}$ lying on the line passing through the origin with slope $\alpha$.
\ed

Basically, $\rays(L)$ encodes how the plane is partitioned by the open cones over the $1$-codimensional faces of the Thurston ball, where the norm is linear.

\bexa If $L$ is a $2$-bridge link, then Theorem \ref{thm: Thurston ball} states that $\rays(L)\subset \{0,\pm 1, \infty\}$.
\eexa

\bc\label{cor: base-type} Let $L$ be a $2$-bridge link. The following conditions are equivalent:
\begin{itemize}
    \item[(1)] The link complement $M_L$ fibers over $S^1$ with fiber $S_{1,0}(L)$;
    \item[(2)] $\rays(L)\subset \{\pm 1\}$;
    \item[(3)] $L$ is base-type.
\end{itemize}
\ec

\begin{figure}[ht]
    \centering
    \includegraphics[width=0.5\textwidth]{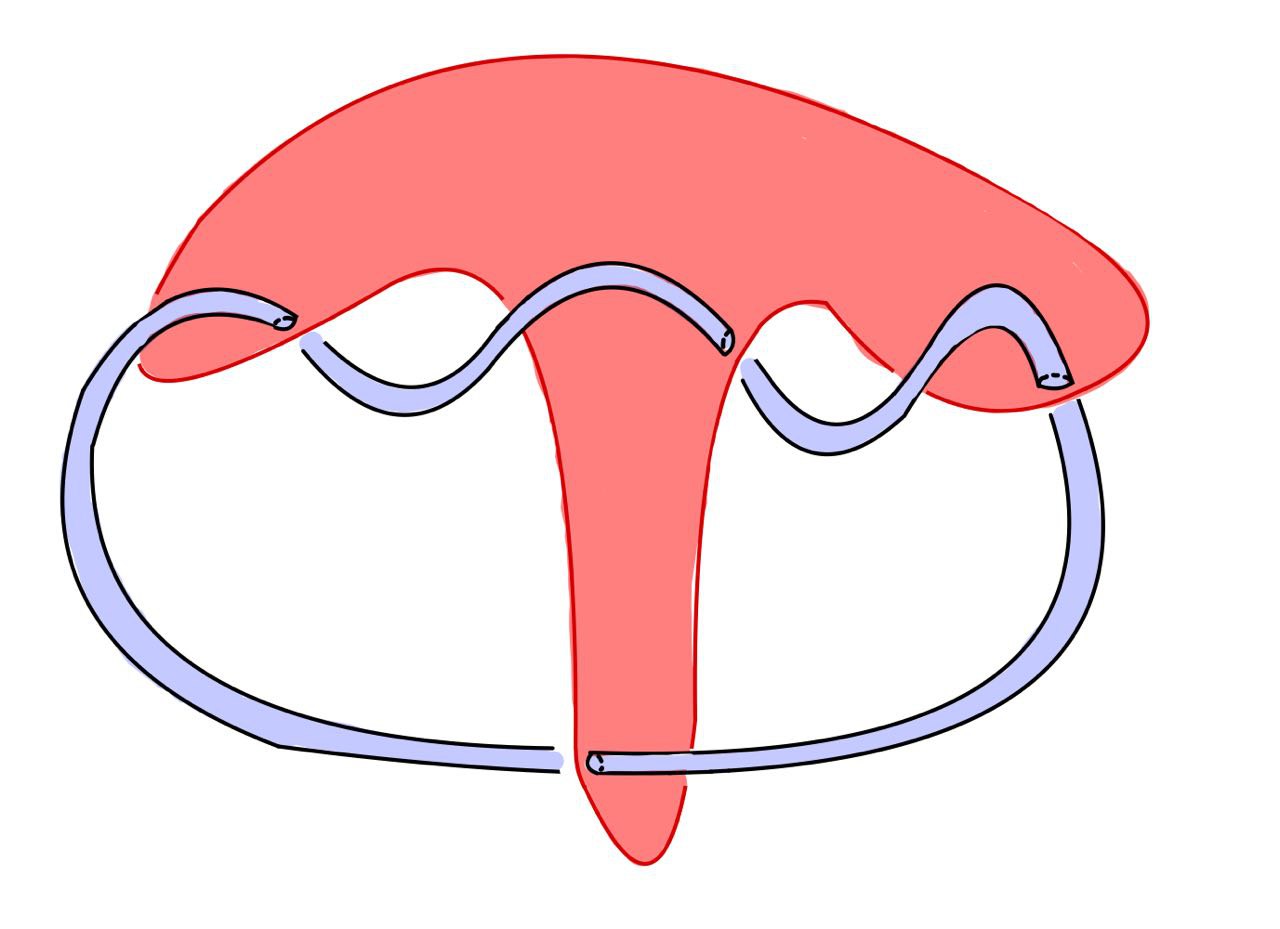}
    \caption{Red: the surface $S_{1,0}(L)$, obtained by removing some discs from a disc $D$ spanning $\ell_1$ in $S^3$. Light blue: the torus $\partial N(\ell_2)$. In the picture, the torus $\partial N(\ell_1)$ is not shown for the sake of clearness.}\label{fig: sut manifold}
\end{figure} 

\bp The implication \emph{(1)$\implies$(2)} follows directly from Theorem \ref{thm: Thurston ball}, the fact that a fibering class must lie in the open cone over a top-dimensional face of the Thurston ball (\cite{norm}), and the symmetries of $L$.

We show the implication \emph{(2)$\implies$(3)} by contradiction. Suppose that $L$ is not base-type. Then, given an alternating diagram for $L$, we can decompose $L$ as a tangle sum of two $2$-bridge links $L_0$ and $L_1$, as in Remark \ref{rem: tangle}.  For $i=0,1$, by abuse of notation, we will call $\ell_1$ and $\ell_2$ the components of $L_i$ coming from $\ell_1$ and $\ell_2$ of $L$ respectively, both with the induced orientations. By equation (\ref{eq: tangel sum}), we get \begin{align*}
    &x_L(\ell_1+\ell_2)+x_L(\ell_1-\ell_2)=x_{L_0}(\ell_1+\ell_2)+x_{L_1}(\ell_1+\ell_2)+2+x_{L_0}(\ell_1-\ell_2)+x_{L_1}(\ell_1-\ell_2)+2=\\ 
    &=\left(x_{L_0}(\ell_1+\ell_2)+x_{L_0}(\ell_1-\ell_2)\right)+ \left(x_{L_1}(\ell_1+\ell_2)+x_{L_1}(\ell_1-\ell_2)\right)+4\ge \\ &\ge 2x_{L_0}(\ell_1)+2x_{L_1}(\ell_1)+4=2x_L(\ell_1)+2>2x_L(\ell_1).
\end{align*}
So, if \emph{(2)} holds, an alternating rational diagram for $L$ cannot have self-crossing boxes, i.e. $L$ must be base-type.

The implication \emph{(3)$\implies$(1)} can be shown via sutured manifold hierarchies. Fix an alternating rational diagram for the base-type link $L$. Start with the sutured manifold $(M,\gamma)=(M_L,\partial M_L)$ and decompose it along the surface $S_{1,0}(L)$, which is a planar surface obtained by removing some discs from a disc $D\subset S^3$ spanning $\ell_1$. See Figure \ref{fig: sut manifold}.

\begin{figure}[ht]
    \centering
    \includegraphics[width=0.5\textwidth]{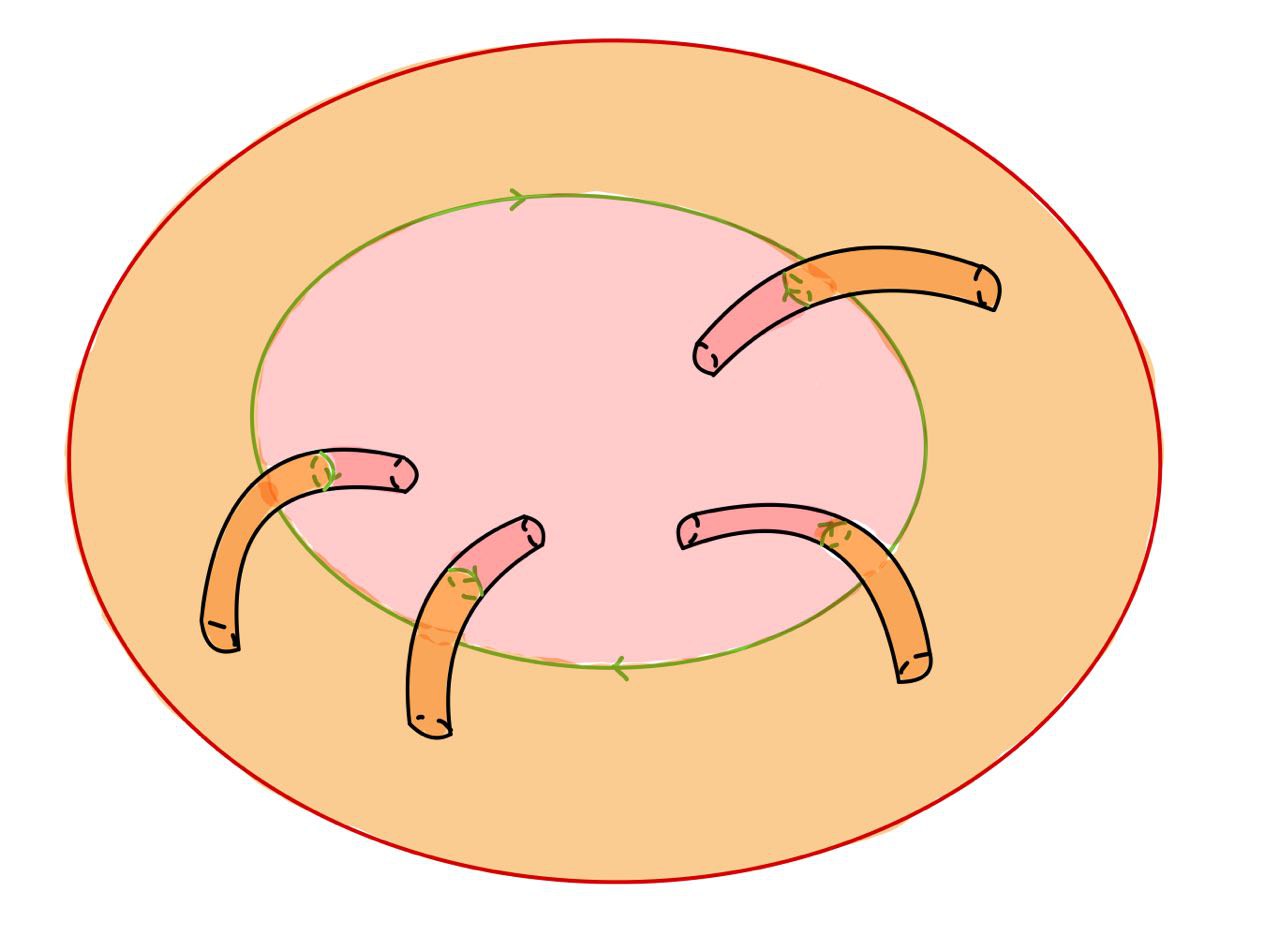}
    \caption{The sutured manifold $(S^3-N,\gamma')$. Green: the sutures $\gamma'$. Pink and orange: the regions $R_+(\gamma')$ and $R_-(\gamma')$.}\label{fig: disc decomp}
\end{figure}

The new sutured manifold is $(S^3-N, \gamma')$, where $N$ is a handlebody given by the union of the ball $B:=D\times [0,1]$ and one solid cylinder $D_i\times [0,1]$ for each arc of $\ell_2$ in the fixed rational diagram of $L$. Each cylinder $D_i\times [0,1]$ is glued to $B$ via maps $D_i\times\{0,1\}\to \partial B$. The equator $\partial D\times {\frac 12}$ is a suture of $\gamma'$, and divides the surface $\partial B- \cup_i D_i\times\{0,1\}$ into two subsurfaces, one belonging to $R_+(\gamma')$ and the other to $R_-(\gamma')$. Since $\ell_2$ intersects $D$ always with the same sign, the gluing curves $\partial D_i\times\{0\}$ and $\partial D_i\times\{1\}$ lie one in $R_+(\gamma')$ and the other in $R_-(\gamma')$. In particular, there is exactly one suture of $\gamma'$ on every tube $\partial D_i\times [0,1]$. After an isotopy, we can put $(S^3-N,\gamma')$ in the form of Figure \ref{fig: disc decomp}. Here, it's easy to see that each handle spans a product-disc in $S^3-N$. Once $(S^3-N,\gamma')$ is decomposed through these product-discs, we obtain the sutured ball $(D^2\times [0,1],\partial D^2\times \{\frac 12\})$. Finally, Theorem \ref{thm: fibering} gives the result.
\ep

A statement analogous to Corollary \ref{cor: base-type} would be that the class $\ell_1+\ell_2$ lies in the open cone over a top-dimensional face of the Thurston ball if and only if $M_L$ fibers over $S^1$ with fiber $S_{1,1}(L)$. Whilst the ``if" part is clearly true, the ``only if" is generally false, see Example \ref{exa: nonfibering}. Nonetheless, knowing that $\ell_1+\ell_2$ lies in the open cone over a top-dimensional face of the Thurston ball gives us precise information about the decomposition of an alternating rational diagram for $L$ into base-type pieces. These pieces will be fairly simple. 

\br Let $L$ be a tangle sum of $L_0$ and $L_1$ as in Remark \ref{rem: tangle}. By equality (\ref{eq: tangel sum}), we have \begin{align*}
    & x_L(\ell_1)+x_L(\ell_2)=x_{L_0}(\ell_1)+x_{L_1}(\ell_1)+1+x_{L_0}(\ell_2)+x_{L_1}(\ell_2)+1=\\ &=\left(x_{L_0}(\ell_1)+x_{L_0}(\ell_2)\right)+\left(x_{L_1}(\ell_1)+x_{L_1}(\ell_2)\right)+2\ge x_{L_0}(\ell_1+\ell_2)+x_{L_1}(\ell_1+ \ell_2)+2=x_L(\ell_1+\ell_2).
\end{align*}

In particular, $x_L(\ell_1)+x_L(\ell_2)=x_L(\ell_1+\ell_2)$ holds if and only if $x_{L_0}(\ell_1)+x_{L_0}(\ell_2)=x_{L_0}(\ell_1+\ell_2)$ and $x_{L_1}(\ell_1)+x_{L_1}(\ell_2)=x_{L_1}(\ell_1+\ell_2)$ both hold.

The analogous statement for $x_L(\ell_1-\ell_2)$ is true as well.
\er

Thanks to the previous remark, our attention is now addressed to base-type $2$-bridge links satisfying either $x_L(\ell_1)+x_L(\ell_2)=x_L(\ell_1+\ell_2)$ or $x_L(\ell_1)+x_L(\ell_2)=x_L(\ell_1-\ell_2)$.

\bl Let $D$ be a base-type diagram for a $2$-bridge link $L$, standardly oriented with components $\ell_1$ and $\ell_2$. Up to mirroring, suppose that every crossing box in $D$ is positive. In $H_2(M_L,\partial M_L)$
\begin{itemize}
    \item $x_L(\ell_1)+x_L(\ell_2)=x_L(\ell_1+\ell_2)$ if and only if $D=T(2), T(1,2k-1)$ or $T(1,2k,1)$ for some $k\in \N$. In this case, we say that $D$ is $(+)$-type.
    \item $x_L(\ell_1)+x_L(\ell_2)=x_L(\ell_1-\ell_2)$ if and only if $D=T(2k)$ or $T(2k-1,1)$ for some $k\in \N$. In this case, we say that $D$ is $(-)$-type.
    \item Both equalities hold if and only if $D=T(2)$ or $D=(1,1)$.    
\end{itemize}
\el
\bp Because of Corollary \ref{cor: base type}, the equality $x_L(\ell_1)+x_L(\ell_2)=x_L(\ell_1+\ell_2)$ holds if and only if $x(\ell_1-\ell_2)=0$. So the lemma follows easily from Proposition \ref{prop: base type 11}. The other case is analogous.
\ep

The next corollary summarises the discussion above, underlying when a $2$-bridge link has hexagonal Thurston ball.

\bc Let $D$ be an alternating rational diagram for a $2$-bridge link $L$, standardly oriented with components $\ell_1$ and $\ell_2$. Up to mirroring, suppose that every crossing box in $D$ is positive. Let $b_1,...,b_r\in \N_+$ be the labels of the self-crossing boxes in $D$, and let $D_i$ be the base-type $2$-bridge link diagram obtained by cutting $D$ along two spheres passing through the boxes $b_i$ and $b_{i+1}$ ($b_0=b_{r+1}=0$), as in Remark \ref{rem: tangle}.
In $H_2(M_L,\partial M_L)$
\begin{itemize}
    \item $\rays(L)\subset \{0,-1,\infty\}$ if and only if, for each $i=0,...,r$, the diagram $D_i$ is $(-1)^{b_0+...+b_i}$-type.
    \item $\rays(L)\subset \{0,1,\infty\}$ if and only if, for each $i=0,...,r$, the diagram $D_i$ is $(-1)^{1+b_0+...+b_i}$-type.   
\end{itemize}
\ec

\bc Let $D$ be a base-type diagram for a $2$-bridge link $L$, standardly oriented with components $\ell_1$ and $\ell_2$. Up to mirroring, suppose that every crossing box in $D$ is positive. Then, $\rays(L)\subset \{0,\infty\}$ if and only if $D$ decomposes as a tangle-sum of $T(2)$ or $T(1,1)$ diagrams.
\ec

\bexa\label{exa: nonfibering} The $2$-bridge link $L=L_{7/16}$ admits the alternating rational diagram $T(2,3,2)$, which is a tangle-sum of two copies of $T(2)$. So, thanks to the previous corollary (or by direct computations), the Thurston ball of $H_2(M_L,\partial M_L)$ has vertices along the axes. In particular, the class $\ell_1-\ell_2$ lies in the open cone over a top-dimensional face of the Thurston ball. 

At the same time, the unique rational diagram of $L$ with even coefficients is $T(2,4,-2)$, so the class $\ell_1-\ell_2$ is not fibered by \cite{gk}. Thus, no analogue of Corollary \ref{cor: base-type} holds in this case.
\eexa

\section{Satellites of $2$-bridge links}\label{sec: satellites}
As an application of the description of the Thurston-norm of $2$-bridge links, we use an iterated satellite construction to yield $2$-component links whose exteriors have arbitrarily complicated Thurston ball.

\bigskip

\begin{figure}[ht]
    \centering
    \includegraphics[width=0.8\textwidth]{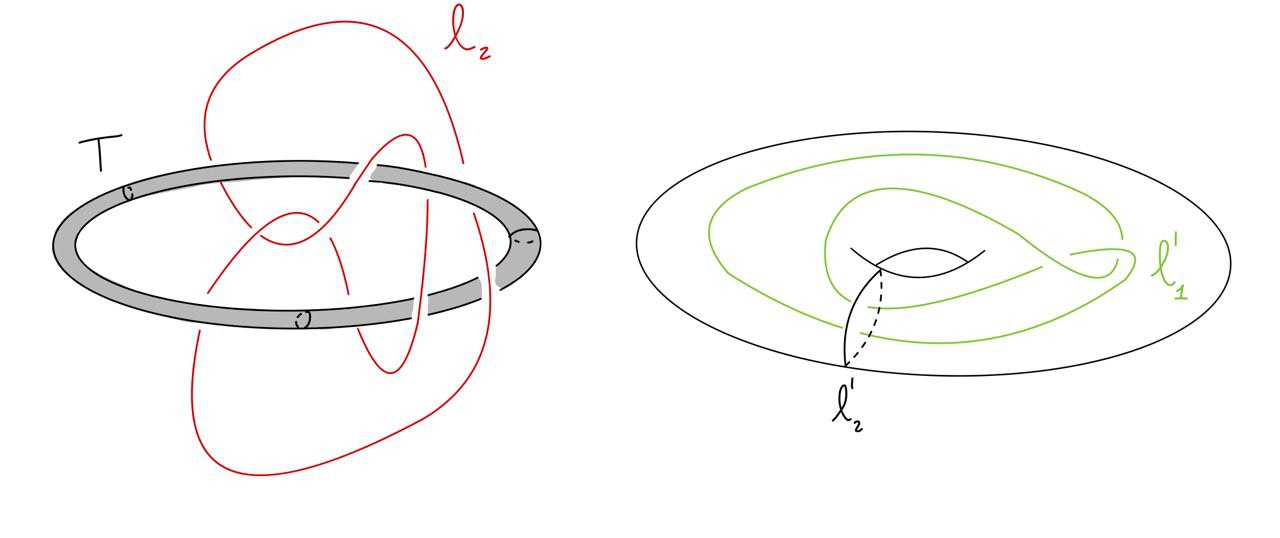}
    \caption{}\label{fig: L and L'}
\end{figure}

Let $L\subset S^3$ be an oriented nonsplit link with components $\ell_1$ and $\ell_2$. Let $\ell_1'\subset D^2\times S^1\subset S^3$ be an oriented knot, and let $L'=\ell_1'\cup\ell_2'$ be the link in $S^3$ with $\ell_2'=\partial D^2\times\{*\}$ an oriented meridional circle. We ask $L'$ to be nonsplit, i.e. $\ell_1'$ is not contained in a $3$-ball of $D^2\times S^1$. See Figure \ref{fig: L and L'}. 

\begin{figure}[ht]
    \centering
    \includegraphics[width=0.5\textwidth]{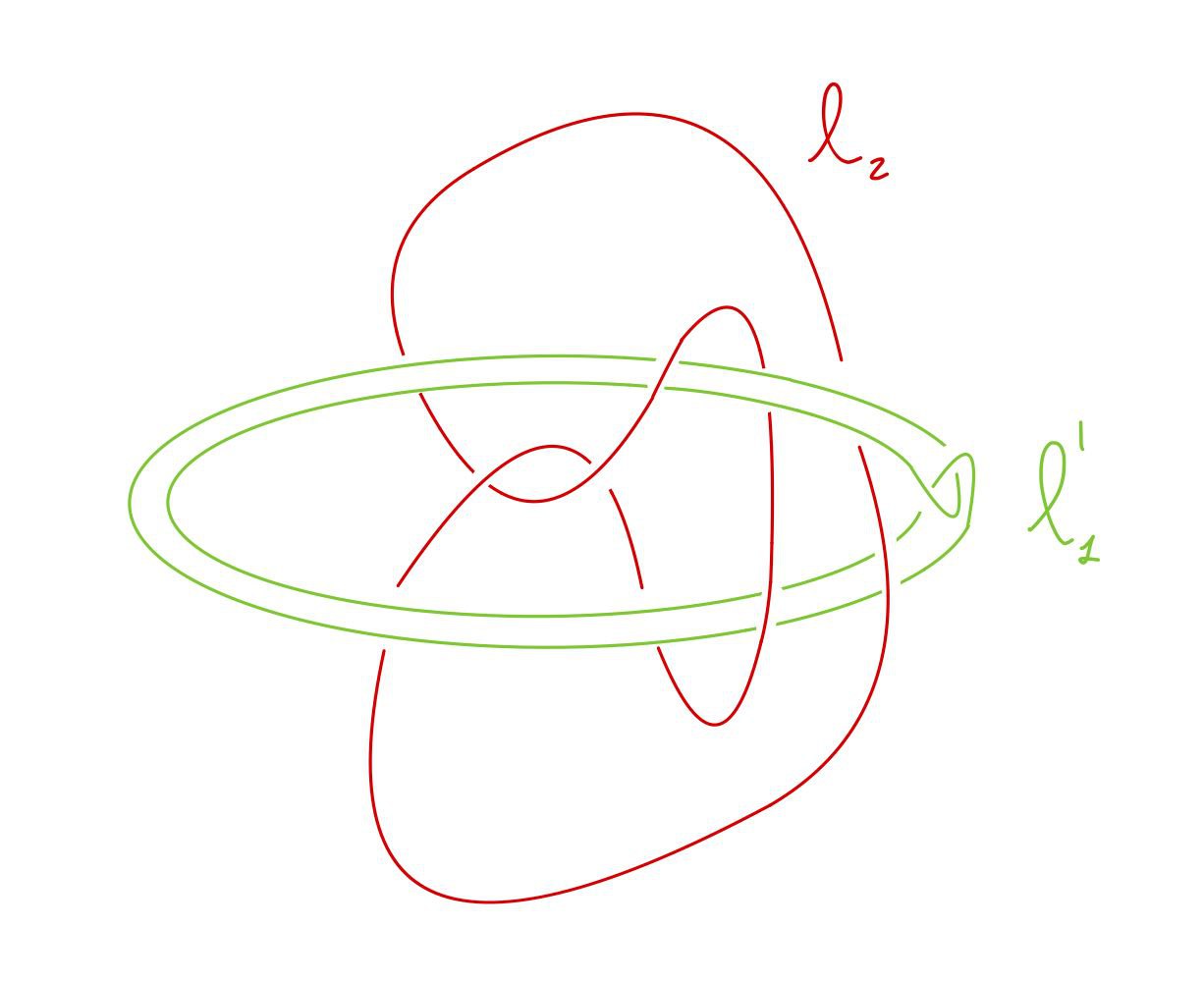}
    \caption{The satellite link $\overline{L}$.}\label{fig: Lbar}
\end{figure}

Let $\lambda$ be the linking number of $\ell_1$ and $\ell_2$, and $\lambda'$ the linking number of $\ell_1'$ and $\ell_2'$. Consider the satellite link $\overline L=\ell_1'\cup \ell_2\subset S^3$ obtained as follows. Fill the boundary component $T=\partial N(\ell_1)\subset M_L$ with the solid torus $D^2\times S^1$ containing $\ell_1'$ so that: \begin{itemize}
    \item $\ell_2'$ is identified with a canonical meridian of $\ell_1$ (i.e. oriented to have linking number $1$ with it) and
    \item a longitude $\{*\}\times S^1\subset \partial D^2\times S^1$ having algebraic intersection $1$ with $\ell_2'$ is identified to a canonical longitude for $\ell_1$. See Figure \ref{fig: Lbar}.
\end{itemize} 
Observe that, in this way, the torus $T$ splits the manifold $M_{\overline L}$ into two pieces $M_L$ and $M_{L'}$, one of which has the reversed orientation, say $M_{L'}$.

\bt\label{thm: satellite norm} For any $a,b\in \R$, the following equality holds $$x_{\overline L}(a,b)=x_L(\lambda' a,b)+x_{L'}(a,\lambda b),$$ where, given an oriented link $P$ with components $p_1$ and $p_2$, $x_P$ indicates the Thurston norm on $H_2(M_P,\partial M_P;\R)$, with basis $p_1,p_2$.
\et
\bp Let $\overline{S}$ be a norm-minimizing surface representing the class $(a,b)=a\ell_1'+b\ell_2\in H_2(M_{\overline L},\partial M_{\overline L})$. We can suppose $\overline S$ is incompressible and has coherently-oriented boundary. After possibly surgeries that do not increase the Euler characteristic, we can also assume that $\overline S$ intersects the torus $T=\partial N(\ell_1)$ in nontrivial parallel closed curves representing the same homology class in $H_1(T;\Z)$. The torus $T$ decomposes $\overline S$ into two surfaces $S=\overline{S}\cap M_L$ and $S'=\overline S\cap M_{L'}$, in such a way that $$\chi(\overline S)=\chi(S)+\chi(S').$$
Since $L$ and $L'$ are nonsplit and the manifolds $M_L$ and $M_{L'}$ are irreducible, we can also ask that no component of $S$ or $S'$ is a sphere or a disk. In particular, the equality $$\chi_-(\overline S)=\chi_-(S)+\chi_-(S')$$ holds. 
The surfaces $S$ and $S'$ must also be norm-minimizing in their homology groups. Indeed, the class $[S]\in H_2(M_L,\partial M_L)$ is uniquely determined as the Poincar\'e dual of the restriction $i^*(PD([\overline S])$, where $i^*:H^1(M_{\overline L})\to H^1(M_L)$ is the map induced by the inclusion. Thus, a coherently-oriented at the boundary representative for $[S]$ with lower $\chi_-$ would have intersection with $T$ isotopic to the one of $S$, so it would patch together with $S'$, then giving a representative for $[\overline S]$ of lower $\chi_-$. By applying the same reasoning to $S'$, we obtain $$x_{\overline{L}}([\overline S])=x_L([S])+x_{L'}([S']).$$ 
To conclude, we just need to express the classes $[S]$ and $[S']$ in terms of the bases for the corresponding homology groups. 

Given an oriented $2$-component link $P=p_1\cup p_2$, the map $\partial: H_2(M_P,\partial M_P)\to H_1(\partial M_P)$ is linear and sends $(1,0)=p_1$ to the direct sum of the homology class (of the canonical longitude for) $p_1$ with $-\lk (p_1,p_2)$ times a canonical meridian for $p_2$. 

Let $b'\in\Z$ be such that $[S']=a\ell_1'+b'\ell_2'=(a,b')$, then $$[S'\cap \partial N(\ell_2')]=b'\ell_2'-\lambda'a\mu_2',$$ where $\mu_2'$ is a canonical meridian for $\ell_2'$. In a similar fashion, let $a'\in \Z$ be such that $[S]=(a',b)$, then $$[S\cap \partial N(\ell_1)]=a'\ell_1-\lambda b\mu_1,$$ where $\mu_1$ is a canonical meridian for $\ell_1$. In the torus $T=\partial N(\ell_1)=-\partial N(\ell_2')$ (where the minus indicates the reversed orientation), the classes $\ell_1$ and $-\mu_2'$ are identified, as well as $\mu_1$ and $-\ell_2'$. So, we get $[S]=(\lambda'a,b)$ and $[S']=(a,\lambda'b)$.
\ep 

\br In the notation above, $$[S'\cap N(\ell_1')]=[\overline S\cap N(\ell_1')]=a\ell_1'-b\lk(\ell_1',\ell_2)\mu_1'.$$ At the same time, $$[S'\cap N(\ell_1')]=a\ell_1'-\lambda'b'\mu_1'=a\ell_1'-\lambda'\lambda b\mu_1'.$$ Thus, the equality $$\lk(\ell_1',\ell_2)=\lambda\lambda'$$ holds.
\er

\bc If $\rays(L)=\{\alpha_1,...,\alpha_k\}$ and $\rays(L')=\{\beta_1,...,\beta_h\}$, then $$\rays(\overline L)=\{\frac {\alpha_1}{\lambda'},...,\frac {\alpha_k}{\lambda'}, \lambda\beta_1,...,\lambda\beta_h\}.$$
\ec
\bp Consider $(a,b),(c,d)\in H_2(M_{\overline L},\partial M_{\overline L})$. The vectors $(a,b)$ and $(c,d)$ lie on the same cone over a face of $x_{\overline L}$ if and only if $$x_{\overline L}(a+c,b+d)=x_{\overline L}(a,b)+x_{\overline L}(c,d).$$ At the same time, by Theorem \ref{thm: satellite norm}, we have \begin{align*}
    &x_{\overline L}(a+c,b+d)=x_L(\lambda' a+\lambda'c,b+d)+x_{L'}(a+c,\lambda b+\lambda d)\le\\ \le & x_L(\lambda' a,b)+x_L(\lambda' c,d)+x_{L'}(a,\lambda b)+x_{L'}(c,\lambda d)=x_{\overline L}(a,b)+x_{\overline L}(c,d)
\end{align*} where equality holds if and only if $(\lambda' a,b)$ and $(\lambda' c,d)$ lie on the same cone over a face for $x_L$, and $(a,\lambda b)$ and $(c,\lambda d)$ lie on the same cone over a face for $x_{L'}$.

For every $\eta\in\R$, consider the linear transformations $T_\eta,T'_\eta:\R^2\to \R^2$ given by $T_\eta(a,b)=(\eta a, b)$ and $T'(a,b)=(a,\eta b)$. 

The discussion above shows that given a cone $C$ for $x_L$ and a cone $C'$ for $x_{L'}$, the intersection of preimages $T^{-1}_{\lambda'}(C)\cap (T')^{-1}_{\lambda}(C')$ is a cone for $x_{\overline L}$. Viceversa, any cone for $x_{\overline L}$ is such an intersection. The corollary follows.
\ep

\bexa If $L$ and $L'$ both have linking number $0$, then the Thurston ball of $\overline L$ has always at most four vertices, lying on the axes.
\eexa

\bexa\label{exa: complicated ball} Let $L$ be a $2$-bridge link with linking number $\lambda\ge 2$ and $\rays(L)=\{0, 1,\infty\}$. Since $L$ has trivial components, we can choose $L'=L=:L_0$ and get $L_1:=\overline L$ with $$\rays(L_1)=\{0, \lambda^{-1},\lambda, \infty\}.$$
We can inductively define a new link $L_i$ by replacing $L\leftarrow L_{i-1}$ and $L'\leftarrow L_0$, obtaining a new $\overline L=:L_i$ with linking number $\lambda^{i+1}$. By induction, we can easily show that $$\rays(L_i)=\{0, \lambda^{-i}, \lambda^{-i+2},..., \lambda^{i-2}, \lambda^i, \infty\}.$$ In particular, $L_n$ is a link with unknotted components whose exterior Thurston ball has $2(n+3)$ faces. Hence, we can state the following. 
\eexa

\bt\label{thm: complexity} For every $n\in\N$, there is a $2$-component link $L\subset S^3$ with each component being the unknot, such that the Thurston ball of $H_2(M_L,\partial M_L;\R)$ has exactly $2n$ faces.
\et

\end{document}